\documentclass[review,10pt]{elsarticle}
\usepackage[utf8]{inputenc}
\usepackage{microtype}
\usepackage[margin=2.5cm]{geometry}
\usepackage{setspace}
\usepackage[none]{hyphenat}
\usepackage[pagewise]{lineno}
\usepackage[english]{babel}
\usepackage{color}
\usepackage{booktabs}
\usepackage{longtable}
\usepackage[most]{tcolorbox}
\usepackage{ragged2e}
\usepackage[overload]{empheq}
\usepackage[justification=justified]{caption}
\usepackage{enumitem}

\usepackage{amsmath,amsfonts,amssymb,mathtools}
\usepackage{cases}
\usepackage{multirow}
\usepackage{graphicx,float}
\usepackage{algpseudocode}
\usepackage{subcaption}
\usepackage{eurosym}
\usepackage{amstext} 
\usepackage{arydshln}
\usepackage{flushend}
\usepackage{threeparttable}
\usepackage{url}

\usepackage[]{hyperref}
\hypersetup{colorlinks, linkcolor=blue, citecolor=blue, urlcolor=blue}

\usepackage[ruled,vlined,linesnumbered]{algorithm2e}
\DeclareRobustCommand{\officialeuro}{%
  \ifmmode\expandafter\text\fi
  {\fontencoding{U}\fontfamily{eurosym}\selectfont e}}
\usepackage{framed} 
\immediate\write18{%
makeindex Main_File.nlo -s nomencl.ist -o Main_File.nls 
}
\usepackage{multicol} 
\usepackage{nomencl}[noprefix] 
\makenomenclature
\renewcommand*\nompreamble{\begin{multicols}{2}}
\renewcommand*\nompostamble{\end{multicols}}

\def\tsc#1{\csdef{#1}{\textsc{\lowercase{#1}}\xspace}}
\tsc{WGM}
\tsc{QE}
\tsc{EP}
\tsc{PMS}
\tsc{BEC}
\tsc{DE}

\begin{document}

    \justifying
    \begin{frontmatter}   
        \title{Adaptive Robust Optimization Models for DER Planning in Distribution Networks under Long- and Short-Term Uncertainties}
    
        \author[1,2]{Fernando García Muñoz\texorpdfstring{\corref{mycorrespondingauthor}}{}} \cortext[mycorrespondingauthor]{Corresponding author}{\ead{fernando.garciam@usach.cl}}{}

        \author[1,2]{Cristian Duran-Mateluna}

        \address[1]{University of Santiago of Chile (USACH), Faculty of Engineering, Industrial Engineering Department, Chile}
        \address[2]{University of Santiago of Chile (USACH), Faculty of Engineering, Program for the Development of Sustainable Production Systems (PDSPS), Chile}

    \begin{abstract}
        This study introduces adaptive robust optimization~(ARO) and adaptive robust stochastic optimization~(ARSO) approaches to address long- and short-term uncertainties in the optimal sizing and placement of distributed energy resources in distribution networks. ARO models uncertainty using a Budget of Uncertainty (BoU), while ARSO distinguishes long-term (LT) demand (via BoU) and short-term (ST) photovoltaics generation (via scenarios). Adapted Benders cutting plane algorithms are presented to tackle the tri-level optimization challenges. The experiments consider a modified version of the IEEE 33 bus system to test these two approaches and also compare them with traditional robust and stochastic optimization models. The results indicate that distinguishing between LT and ST uncertainties using a hybrid formulation such ARSO yields a solution closer to the optimal solution under perfect information than ARO.
    \end{abstract}


\end{frontmatter}

\section{Introduction}\label{Intro_Chap}

The integration of Distributed Energy Resources (DERs) into distribution networks (DNs) has multi-faceted advantages, such as bolstering power system resilience, improving energy efficiency, and reducing reliance on conventional energy sources~\cite{bollen2011integration}. DERs include distributed generation (DG), such as photovoltaic (PV) generation and wind turbines, as well as battery energy storage systems (BESS), all of which require efficient planning to optimize their deployment from both a technical and economic perspective~\cite{9698197}. However, one of the most significant challenges in DN planning and operation is balancing electricity demand with renewable energy generation \cite{9233134}. Thus, the integration of wind and PV systems into transmission and DNs has driven significant research efforts to manage their inherent uncertainties and ensure system balance~\cite{7937825}.

In this context, effectively integrating DERs into DNs requires a comprehensive planning approach that considers infrastructure design, the deployment of renewable energy resources, and the management of uncertainties in power system operations~\cite{9233331}. As power systems grow in complexity, the adoption of optimization models capable of addressing these uncertainties has become essential for strategic decision-making \cite{milligan2012stochastic}. In this regard, \emph{stochastic optimization} has emerged as a valuable tool for managing uncertainty based on the probability distribution of random variables of uncertain parameters. When the probability distribution is well-defined, solutions with high accuracy can be obtained by considering a large number of scenarios, albeit at the cost of significant computational effort~\cite{de2014computational}. Alternatively, when information about the distribution of probabilities is unavailable or sparse, \emph{robust optimization} uses predefined uncertainty sets to ensure feasible solutions under worst-case conditions,  reducing computational complexity but often leading to conservative decisions that may impact investment decisions~\cite{gorissen2015practical}. This trade-off between computational efficiency and solution conservatism has been studied and applied to tackle the challenge of optimally integrating renewable energy sources into power systems~\cite{lorca2016robust}, balancing the need for reliability, feasibility, and cost-effectiveness. In both approaches exist single-stage, two-stage, or multi-stage formulations, which are defined according to the point in time at which decisions are made regarding the revealing of the uncertainty~\cite{Sun2021}. 

Despite the existence of different formulations to handle uncertainty, its behavior can significantly vary over different planning horizons~\cite{9739994}. For instance, uncertainty related to PV generation is characterized by intraday variability, often constrained to a 24-hour horizon, while long-term (LT) electricity demand could be more challenging to estimate due to the adoption of emerging technologies, such as electric vehicles and the increasing electrification of other energy sectors~\cite{BARINGO2020114679}. In this context, it becomes essential to distinguish between LT uncertainties, such as annual electricity demand projections, and short-term (ST) uncertainties, such as PV power generation variability~\cite{zhang2018adaptive,7954668}. Addressing these uncertainties effectively in planning models is crucial for supporting decision-makers in both investment and operational cost assessments, ensuring network reliability, and accommodating the expected growth in electricity demand over the coming years~\cite{REN2024110378}.

To the best of our knowledge, unlike classical network expansion problems, the optimal sizing and placement of DERs (OSP-DER) in DN under uncertainty conditions, considering both LT and ST uncertainties, has not been extensively studied. The contributions of this paper are threefold: (1) the introduction of an \emph{adaptive robust optimization~(ARO)} and an \emph{adaptive robust stochastic optimization~(ARSO)} formulations for the OSP-DER of PV systems and BESS in DNs, explicitly distinguishing between LT uncertainties in electricity demand and ST uncertainties in PV generation; (2) the adaptation of the Benders decomposition algorithms to efficiently manage the budgeted uncertainty sets in order to solve the ARO and ARSO dual subproblems; and~(3) a comprehensive comparative analysis between stochastic and robust optimization approaches, including \emph{two-stage stochastic optimization}, \emph{single-stage robust optimization}, ARO, and ARSO formulations. 


The remainder of the paper is organized as follows. Section~\ref{sec:literature} presents the literature review on planning problems under uncertainty. Section~\ref{sec:deterProblem} describes the deterministic problem and its two-stage stochastic formulation. Section~\ref{sec:robustProblem} introduces its robust reformulations and the budget of uncertainty sets. Section~\ref{sec:algo} explains the adapted Benders decomposition algorithms for our robust models. Finally, Section~\ref{sec:results} discusses the case study and computational results, while Section~\ref{sec:conclusions}  concludes with the main findings and future research.

\section{Literature Review}\label{sec:literature}
 
Among the different optimization approaches under uncertainty, two-stage stochastic optimization (TSSO) has been widely applied in designing and operating energy systems to manage uncertainties and optimize decision-making processes effectively. For instance,~\cite{han2021two} introduces a two-stage approach for multi-microgrid systems, focusing on minimizing both capital costs for cable installation and operational expenses. Similarly,~\cite{ZHENG2021110835} adopts this method to handle the load and PV production uncertainties for optimal residential PV and BESS systems, where the first stage optimizes the investment portfolio and the second stage minimizes operational costs. In the context of DN, \cite{ortiz2018stochastic} employs a TSSO model for expansion planning (EP) under wind and load uncertainties, while~\cite{NAZARI2019281} applies this approach to conduct a cost-benefit analysis of BESS installations in microgrids (MG). Additionally, \cite{BOFFINO2019104457} uses TSSO to determine the optimal mix of generation and transmission capacities, minimizing investment costs for future demands. Beyond two-stage models, \cite{9234713} extends the concept to multi-stage stochastic programming to address DG uncertainties over multiple periods, introducing a nested decomposition method to manage computational complexity. Furthermore, \cite{GARCIAMUNOZ2022100624} utilizes a TSSO model for the OSP-DER problem, considering electric vehicles and demand response.

On the other hand, in an early stage, some studies explored the use of {single-robust optimization (SRO)} as an approach to address uncertainty, particularly in cases where the probability distribution of random variables is difficult to define. For instance, \cite{6547161,6784363} apply SRO to manage renewable generation and load uncertainty in transmission network EP problems. However, SRO often produces conservative solutions, which may lead to overly cautious decisions in planning problems. To mitigate this limitation, the literature shifted its focus to explore the Adaptive Robust Optimization (ARO)~\cite{7098436}. ARO introduces the second stage to enable decisions to adapt dynamically to the worst-case scenario, providing a more flexible and less conservative approach to uncertainty management. For example, ARO has been used to address uncertainty in operational challenges such as unit commitment~\cite{lee2013modeling,8636235} and economic dispatch~\cite{8293834,zhang2016robust}, as well as in planning problems like EP~\cite{xie2020two}.
In particular, the work in~\cite{ELMELIGY2024110164} presents an ARO model for transmission EP, integrating network-constrained unit commitment as a third-level decision and solving the resulting three-level problem using a column-and-constraint generation (CCG) algorithm. Similarly, the study in~\cite{RADMANESH2024110062} formulates a three-stage, four-level min–min–max–min problem to plan new transmission lines and identify profitable locations for merchant energy storage investments, leveraging a modified CCG to manage the nested structure.

In the context of DN,~\cite{BARINGO2020114679} proposes an ARO-based EP model to optimize the construction of renewable generation units, storage facilities, and charging stations for electric vehicles. \cite{BARINGO2020114679} accounts for ST variability in demand and LT uncertainty in electricity prices, reflecting a nuanced approach to uncertainty management. Similarly, \cite{zhang2018adaptive} employs an LT and ST uncertainty distinction to address the transmission EP problem. Likewise, the study in~\cite{8012502} introduces an ARO EP model for DN, aiming to determine the timing of feeder reinforcements as well as the location, capacity, and installation schedule of dispatchable and wind-based DERs. These efforts highlight the growing sophistication of ARO methods in tackling the inherent complexities of power system expansion under uncertainty.

In previous ARO formulations, uncertainties are often categorized into LT and ST types, reflecting their distinct characteristics and time scales. LT uncertainty typically encompasses high variability and unknown probability distributions over extended periods~\cite{7954668}, such as annual electricity demand projections in DN with high electrification levels. In contrast, ST uncertainty relates to daily operations~\cite{REN2024110378}, such as PV power generation variability. This differentiation is crucial for addressing planning problems effectively, as it allows for more precise modeling of uncertainties while avoiding excessive complexity in the mathematical formulation. Thus, some works have explored hybrid approaches like in~\cite{ALNOWIBET2024110546, 7944676,7954668}, where adaptive robust stochastic optimization (ARSO) approaches for the generation and transmission EP problem are presented considering LT and ST uncertainty. However, these articles are focused on the transmission networks and do not address the OSP-DER in DN.

In relation to solution methods, many metaheuristic methods and specialized algorithms have been developed to address the OSP-DER problem in DN, aiming to explore large solution spaces efficiently and address computational challenges inherent in these problems. For instance,~\cite{RAMADAN2023101872} employs the Artificial Hummingbird Algorithm under uncertainty conditions for optimal DG allocation, while~\cite{8316645} applies the Ant Lion Optimizer to determine the optimal size and placement of DGs. The Artificial Bee Colony Algorithm is utilized in~\cite{6029809} to minimize total system real power loss while solving for DG placement and sizing. Similarly, \cite{kollu2014novel} implements the Harmony Search Algorithm, \cite{SUDABATTULA2016270} employs the Bat Algorithm, and~\cite{WONG2019100892} utilizes the Whale Optimization Algorithm. Other approaches include the Firefly Algorithm in~\cite{anbuchandran2022multi} and the Particle Swarm Optimization Algorithm, which has been extensively applied in studies such as~\cite{10194717,10488786, STAWFEEK201886}. While metaheuristic methods are valuable for addressing computational challenges, the complexity of nested structures of previous robust problems has been driven mainly by Benders decomposition~\cite{jiang2011robust, zhang2018adaptive} and CCG algorithms~\cite{ZENG2013457}.   Benders decomposition method, in particular, is emphasized for its adaptability and efficiency in addressing LT and ST uncertainties, especially for network expansion models~\cite{jiang2011robust}.

While stochastic and robust optimization techniques have been extensively applied to power system planning, such as the EP problem, the challenge of integrating DERs in DNs considering both LT and ST uncertainties remains almost unexplored. Existing studies have either assumed predefined probability distributions or adopted conservative uncertainty sets, often leading to suboptimal investment decisions. To address this gap, this work adapts the ARO and ARSO approaches to distinguish LT and ST uncertainties explicitly, providing a decision-making framework that mitigates over-investment risks. 


\newpage
\section{OSP-DER problem formulation} \label{sec:deterProblem}

In this section, we first present the deterministic OSP-DER problem, followed by the uncertainty parameters with the two-stage stochastic optimization formulation.  The main notation used in this paper is stated below for quick reference, while other symbols are defined as needed throughout the text.

\begin{table}[!h]   
\small
\begin{framed}
\nomenclature[01]{\textbf{Sets}}{}
\nomenclature[02]{$i \in \mathcal{B}$}{Set of buses.}
\nomenclature[04]{$t \in \mathcal{T}$}{Set of time slots.}
\nomenclature[05]{$s \in \mathcal{S}$}{Set of scenarios.}
\nomenclature[06]{$(i,j) \in \mathcal{L}$}{Set of lines, such that $\mathcal{L} = \lbrace (i,j);i,j \in \mathcal{B}\rbrace$.}
\nomenclature[07]{\textbf{Parameters}}{}
\nomenclature[08]{$PG^{max}_{i}$}{Maximum active power output at bus $i$ [kW].}
\nomenclature[08]{$QG^{max}_{i}$}{Maximum reactive power output at bus $i$ [kvar].}
\nomenclature[08]{$QG^{min}_{i}$}{Minimum reactive power output at bus $i$ [kvar].}
\nomenclature[09]{$PL_{i,t,s}$}{Active load of bus $i$ at period $t$ in scenario $s$ [kW].}
\nomenclature[10]{$QL_{i,t,s}$}{Reactive load of bus $i$ at period $t$ in scenario $s$ [kvar].}
\nomenclature[11]{$V^{min}_i$}{Minimum voltage allowed at bus $i$ [p.u.].}
\nomenclature[12]{$V^{max}_i$}{Maximum voltage allowed at bus $i$ [p.u.].}
\nomenclature[14]{$S^{max}_{i,j}$}{Allowed maximum apparent power of the line between buses $i$ and $j$ [kVA].}
\nomenclature[17]{$PV_{t,s}$}{Maximum PV active power available to inject at period $t$  in scenario $s$[\%].}
\nomenclature[18]{$SOC^{min}$}{Minimum state of charge of battery [\%].}
\nomenclature[19]{$SOC^{max}$}{Maximum state of charge of battery [\%].}
\nomenclature[20]{$\varphi^{ch}$}{Charging battery efficiency [\%].}
\nomenclature[21]{$\varphi^{ds}$}{Discharging battery efficiency [\%].}
\nomenclature[22]{$PB$}{Maximum charge/discharge power for battery [kW].}
\nomenclature[24]{$\Delta t$}{Time slot. (1 h in the model).}
\nomenclature[27]{$\lambda_{t,s}$}{Price of the energy bought from the grid at time $t$ in scenario $s$.}
\nomenclature[29]{$\rho_{s}$}{Probability of the scenario $s$.}
\nomenclature[30]{$C^{pv}_{i}$}{Maginal cost of installed capacity in PV at bus $i$.}
\nomenclature[31]{$C^{bt}_{i}$}{Maginal cost of installed capacity in BESS at bus $i$.}
\nomenclature[32]{$I^{pv}_{i}$}{Cost of installing a PV system at bus $i$.}
\nomenclature[33]{$I^{bt}_{i}$}{Cost of installing a BESS at bus $i$.}
\nomenclature[34]{$OC^{pv}$}{PV operational cost.}
\nomenclature[35]{$OC^{bt}$}{BESS operational cost.}
\nomenclature[36]{$N^{pv}$}{Maximum allowable PV systems for installation.}
\nomenclature[37]{$N^{bt}$}{Maximum allowable BESS for installation.}
\nomenclature[38]{$\Upsilon$}{Maximum and minimum capacity for installation for BESS or PV systems.}
\nomenclature[39]{$\beta_{pl,pv}$}{Budget of Uncertainty of power load or PV.}
\nomenclature[40]{$\widetilde{PV}_{t}$}{Uncertainty parameter at period $t$.}
\nomenclature[41]{$\overline{PV}_{t}$}{Expected value of the Uncertainty parameter at period $t$.}
\nomenclature[42]{$\widehat{PV}_{t}$}{Deviation of the expected value at period $t$.}
\nomenclature[44]{\textbf{Variables}}{}
\nomenclature[45]{$\nu^{pv}_{i}$}{Binary variable: 1 if a PV system is installed in bus $i$, 0 otherwise.}
\nomenclature[46]{$\nu^{bt}_{i}$}{Binary variable: 1 if a BESS is installed in bus $i$, 0 otherwise.}
\nomenclature[47]{$\gamma^{pv}_{i}$}{PV capacity installed at bus $i$.}
\nomenclature[48]{$\gamma^{bt}_{i}$}{BESS capacity installed at bus $i$.}
\nomenclature[49]{$\Delta p_{i,t,s}$}{Difference between the energy self-generated and the consumption in bus $i$ at period $t$ in scenario $s$ [kW].}
\nomenclature[50]{$pv_{i,t,s}$}{Power self-generated in bus $i$ at period $t$ in scenario $s$ [kW].}
\nomenclature[50]{$pg_{i,t,s}$}{Power at Substation in bus $i$ at period $t$ in scenario $s$ [kW].}
\nomenclature[51]{$v_{i,t,s}$}{Voltage of bus $i$ at period $t$ in scenario $s$ [p.u].}
\nomenclature[52]{$\theta_{i,t,s}$}{Angle of bus $i$ at period $t$ in scenario $s$.}
\nomenclature[52]{$p^{i,j}_{t,s}$}{Active power in line between buses $i$ and $j$ at period $t$ in scenario $s$ [kW].}
\nomenclature[53]{$q^{i,j}_{t,s}$}{Reactive power in line between buses $i$ and $j$ at period $t$ in scenario $s$ [kvar].}
\nomenclature[54]{$ch_{i,t,s}$}{Power charged by the battery in the bus $i$ at period $t$ in scenario $s$ [kW].}
\nomenclature[55]{$ds_{i,t,s}$}{Power discharged by the battery in the bus $i$ at period $t$ in scenario $s$ [kW].}
\nomenclature[56]{$soc_{i,t,s}$}{State of charge of battery of bus $i$ at period $t$ in scenario $s$ [kWh].}
\nomenclature[57]{$w_{i,t,s}$}{Binary variable: 1 if the battery is charging in bus $i$ at period $t$ in scenario $s$, 0 otherwise.}
\input{Article.nls}
\end{framed}
\end{table}

\newpage
\subsection{Deterministic problem}\label{Det_Model_Subsection}

The deterministic OSP-DER problem can be formulated as a Mixed-Integer Linear Programming (MILP) model. Let us consider a radial DN represented as a directed graph ${\mathcal{G}} = ({\mathcal{B}},{\mathcal{L}})$ where ${\mathcal{B}}$ denotes the set of buses and ${\mathcal{L}}$ represents the set of distribution lines, each connecting a pair of adjacent buses $i$ and $j$. This article addresses the planning problem of the PV systems and BESS within a time horizon represented by a set of time slots $ \mathcal{T}$.  The objective is to reduce dependence on the upstream grid while ensuring that electricity demand is met through grid purchases and local generation at each period of time, with all power exchanges respecting distribution line thermal limits. Consequently, the objective function \eqref{OB_Det_Problem} minimizes the investment ($IC_{i}$) and the operation cost ($OC_{i,t}$) at each bus $i \in \mathcal{B}$ and time slot $t \in \mathcal{T}$.

\begin{flalign}\label{OB_Det_Problem}
&\min \enspace z = 
\displaystyle\sum_{i\in {\mathcal {B}}} IC_{i} + \displaystyle\sum_{i\in {\mathcal {B}}}  \displaystyle\sum_{t\in {\mathcal {T}}}  OC_{i,t}
\end{flalign}

For each bus $i \in \mathcal{B}$, the installation decisions are represented by the binary variables $\nu^{pv}_{i}$ and $\nu^{bt}_{i}$, while their corresponding capacities are represented by the continuous variables $\gamma^{pv}_{i}$ and $\gamma^{bt}_{i}$, where  $pv$ and $bt$ identify the investments in PV systems and BESS, respectively.  In the total investment costs, we differentiate marginal capacity cost  $C_i$ associated  with variables $\gamma_{i}$, and the installation cost $I_i$  associated with $\nu_{i}$ as follows:

\begin{flalign}
    & IC_{i} = \gamma^{pv}_{i} C^{pv}_{i} + \nu^{pv}_{i}I^{pv}_{i} + \gamma^{bt}_{i} C^{bt}_{i} + \nu^{bt}_{i}I^{bt}_{i} \label{IC_DM} 
\end{flalign}

Meanwhile, for the operational costs at each bus $i \in \mathcal{B}$ and time slot $t \in \mathcal{T}$, we consider the power supplied by the main grid $pg_{i,t}$ and its price $\lambda_{i,t}$, the power provided by the PV system $pv_{i,t}$, and the power discharged from the BESS $ds$  with their corresponding operational costs $OC^{pb}$ and $OC^{bt}$, as follows:

 \begin{flalign}
    & OC_{i,t} = \lambda_{i,t} pg_{i,t} + OC^{pv}pv_{i,t} + OC^{bt}ds_{i,t} \label{OC_DM}
\end{flalign}

A feasible solution corresponds to a subset of buses that satisfies network and power generation constraints. For each $i \in {\mathcal {B}}$, the installation constraints for PV systems and BESS are defined as follows:

\begin{subequations}\label{Installation_Constraints}
\begin{flalign} 
&\Upsilon^{min} \nu^{pv}_{i} \leq \gamma^{pv}_{i} \leq \Upsilon^{max} \nu^{pv}_{i} \label{PV_Capacity} \\
&\displaystyle\sum_{i \in {\mathcal{B}}} \nu^{pv}_{i} \leq N^{pv} \label{N_PVs} \\
&\Upsilon^{min} \nu^{bt}_{i} \leq \gamma^{bt}_{i} \leq \Upsilon^{max} \nu^{bt}_{i}  \label{Bat_Cap} \\
&\displaystyle\sum_{i \in {\mathcal{B}}} \nu^{bt}_{i} \leq N^{bt} \label{N_Bat_C} 
\end{flalign}
\end{subequations}

Constraints~\eqref{PV_Capacity} limit the installed capacity of PV systems $\gamma^{pv}_i$ by the minimum and maximum capacity thresholds $\Upsilon$, in relation to the installation decisions $\nu_i^{pv}$, while Constraint~\eqref{N_PVs} limits the total number of PV installations in the network. Similarly, Constraint~\eqref{Bat_Cap} bounds the storage capacity of BESS $\gamma^{bt}_i$ in relation to $\nu_i^{bt}$, while Constraint~\eqref{N_Bat_C} restricts the total number of storage units that can be installed.

The classical representation of network constraints through the AC optimal power flow (AC-OPF) equations is known to be NP-hard~\cite{bienstock2019strong}. To ensure computational tractability in the planning problem, we adopt a second-order cone relaxation~\cite{bansalconvex,molzahn2019survey},  and we introduce an additional assumption that network losses are neglected, based on the premise that these losses remain below 5\%. This simplification is reasonable in DNs where line lengths are typically short and power flows are lower than in high-voltage transmission systems. Moreover, the model assumes that local generation sources are geographically close to consumption nodes, which reduces real power losses by minimizing energy transport across the network and helps maintain line loadings within acceptable limits. This limitation is justifiable in a planning context, where the primary objective is to analyze the impact of uncertainty rather than achieving an accurate representation of network losses. Under these assumptions, the operational network constraints for each node $ i \in {\mathcal {B}}$ and time slot $t \in {\mathcal {T}}$  are defined as follows:

\begin{subequations}\label{Network_Constraints}
\begin{flalign} 
   && &\displaystyle\sum_{(i,j) \in {\mathcal{L}}} p^{i,j}_{t} -  \displaystyle\sum_{(j,i) \in {\mathcal{L}}} p^{j,i}_{t} = \Delta p_{i,t}  \label{Active_Balance}&&  \\
    &&&\displaystyle\sum_{(i,j) \in {\mathcal{L}}} q^{i,j}_{t} - \displaystyle\sum_{(j,i) \in {\mathcal{L}}} q^{j,i}_{t} = \Delta q_{i,t} \label{Reactive_Balance}&&  \\
   && &pg_{i,t} \leq PG^{max}_{i} \label{PG_MAX} &&  \\
    &&&QG^{min}_{i} \leq qg_{i,t} \leq QG^{max}_{i} \label{QG_MINMAX} &&  \\
    &&&V^{min}_{i} \leq v_{i,t} \leq V^{max}_{i} \label{V_Limits} &&  \\
    &&&\Delta p_{i,t} = pg_{i,t} + pv_{i,t} - PL_{i,t} + ds_{i,t} - ch_{i,t} \label{DeltaP}\\
   && &\Delta q_{i,t} = qg_{i,t} - QL_{i,t} \label{DeltaQ}&&  \\
   && &v_{j,t} = v_{i,t} -2 (R_{i,j} p^{i,j}_{t} + X_{i,j}q^{i,j}_{t}) \quad && \forall (i,j) \in \mathcal{L}\label{Line_Balance}   \\   
   && &p^{i,j}_{t}A_{r} + q^{i,j}_{t}B_{r} + S^{max}_{i,j}C_{r} \leq 0 \quad && \forall (i,j) \in \mathcal{L}\label{S_Max_Limit}  \\
   && &\theta_{i,t} -  \theta_{j,t} = X_{i,j} p^{i,j}_{t} - R_{i,j}q^{i,j}_{t} \quad & &\forall (i,j) \in \mathcal{L}\label{Dif_Angle} \\
   && &pv_{i,t} \leq PV_{t} \gamma^{pv}_{i} \label{PV_Injection}&  \\
	&&&soc_{i,t+1}= soc_{i,t}+[\varphi^{ch} ch_{i,t}-\frac{1}{\varphi^{ds}}ds_{i,t}]\Delta t  \label{SOC_C}&&  \\ 
	&&&SOC^{min} \gamma^{bt}_{i} \leq soc_{i,t} \leq SOC^{max} \gamma^{bt}_{i}  \label{SOC_Limits}&&  \\
	&&&ch_{i,t}  \leq PB (w_{i,t})  \label{Ch_c}&&  \\ 
	&&&ds_{i,t} \leq PB (1-w_{i,t}) - PB(1- \nu^{bt}_{i})  \label{Ds_C}&&  \\
   && &w_{i,t} \leq \nu^{bt}_{i}  \label{Charging_Var}
\end{flalign}
\end{subequations}

Constraints~\eqref{Active_Balance} and \eqref{Reactive_Balance} define the active and reactive power balance at each node, excluding losses as proposed in~\cite{baran1989network}. Constraints~\eqref{PG_MAX}-\eqref{DeltaQ} impose limits on active power, reactive power, and voltage levels. Constraints~\eqref{DeltaP} capture the aggregated effect of consumption and self-generation, with a similar formulation applied in Constraints~\eqref{DeltaQ}. Likewise, Constraint~\eqref{Line_Balance} represents the line balance equation without accounting for losses in relation to the active $PL_{i,t}$ and reactive loads $QL_{i,t}$. Following~\cite{yuan2019second}, Constraint~\eqref{S_Max_Limit} provides a linearized expression for the maximum line capacity, while Constraints~\eqref{Dif_Angle} incorporates angle variations into the model. Constraints~\eqref{PV_Injection} limit the power generated by the solar panels $pv_{i,t}$ by the installed solar capacity $\gamma^{pv}_i$ and the normalized solar power generation profile $PV_t$. Similarly, Constraints~\eqref{SOC_C}-\eqref{Ds_C} describe the battery operation in relation to state of charge $soc_{i,t}$, the power charged $ch_{i,t}$, or discharged $ds_{i,t}$ by the battery. Binary variables $w_{i,t}$ are used to avoid simultaneous charging and discharging, which is limited by the installation decisions $\nu_i^{bt}$ in Constraints~\eqref{Charging_Var}.

\subsection{Uncertainty parameters and two-stage stochastic formulation}

The optimal integration of DERs in a DN is subject to multiple uncertainties, such as PV generation variability, electricity demand fluctuations, and future prices. This article exclusively addresses the uncertainty associated with electricity demand and power generation from PV systems, leaving the uncertainties related to electricity prices and technology costs for future research. As mentioned in Section \ref{sec:literature}, the main approach used in literature has been the TSSO. In this framework, first-stage (here-and-now) decisions correspond to investment decisions that must be determined before uncertainty is revealed. Meanwhile, second-stage (wait-and-see) decisions represent operational variables, which can be adjusted once uncertain parameters, such as electricity demand and PV generation, become known. Consequently, in this approach, the objective function minimizes investment costs and the expected value of the operational costs.

To provide a compact TSSO formulation of the OSP-DER, let $x=(\nu_i^{pv}, \nu_i^{bt},  \gamma_i^{pv}, \gamma_i^{bt})$ be the first-stage variables of investments decisions,   $\Phi$ be its corresponding feasibility set defined by Constraints~\eqref{PV_Capacity}-\eqref{N_Bat_C}, \( u = ( \widetilde{PL}_{i,t}, \widetilde{PV}_t^{max} ) \)  the uncertain parameters associated with electrical load variations and the power generated by the PV systems,  \mbox{$y= (pg_{i,t}, pv_{i,t}, ds_{i,t}, ch_{i,t}, v_{i,t}, w_{i,t}, p_{i,t}, q_{i,t}, qg_{i,t}, \theta_{i,t} )$} be the second-stage variables of the operational decisions, which are determined by  Constraints~\eqref{Active_Balance}-\eqref{Charging_Var}, and $ C_{I}^{\top}$ and $C_{O}^{\top}$ be the investment and operational cost coefficient vectors, respectively. In the TSSO formulation, the probability distribution of $u$ is assumed to be known, hence the expected values of the uncertain parameters is minimized in the objective function. Consequently, the  TSSO OSP-DER  can be formulated  as follows:
\vspace{-0.4cm}
\begin{subequations}\label{2S_Compact}
\begin{align}[left={(TSSO)}\empheqlbrace]
 \nonumber \underset{x\in \Phi}{\min} \enspace  C_{I}^{\top}x +  \mathbb{E}_u[Q(x,u)] \\
   \nonumber \text{where} \quad Q(x,u) = \min_{y} \enspace &  C_{O}^{\top} y \\
   s.t. \enspace  &h(x,y,u) = 0 \label{Snd_stage_TSSO_1}\\
    &g(x,y,u) \leq 0 \label{Snd_stage_TSSO_2}\\
   \nonumber  & y \in \mathbb{R} \label{var_Secod_stage_TSSO}
\end{align}
\end{subequations}

\noindent where  Constraints~\eqref{Snd_stage_TSSO_1} and ~\eqref{Snd_stage_TSSO_2} encapsulate the operational constraints related to the second-stage variables. Usually, the probability space is discretized into a finite set of scenarios $s \in S$, representing different realizations of the uncertain parameter~\citep{birge2011introduction}. In the nomenclature, all corresponding second-stage variables and parameters include this index. In the following sections, we present robust reformulations from this compact formulation.

\section{Robust optimization approach}\label{sec:robustProblem}

In robust optimization,  the probability distribution of the uncertain parameters is considered unknown, so the uncertainty is then modeled using uncertainty sets, which define the range of possible values that the uncertain parameters can take. Common types of uncertainty sets include Box, Ellipsoidal, Budgeted, Polyhedral, and Norm-Based sets, each offering different trade-offs between robustness and computational complexity~\cite{Sun2021}. This approach aims to identify the worst-case scenario within the defined uncertainty set, ensuring system feasibility and reliability under all possible uncertain outcomes.  In the following, we present how to reformulate the robust OSP-DER in DNs considering budgeted uncertainty sets.

\subsection{Single-Stage Robust Optimization (SRO) formulation}

In this formulation, all decisions must be made before the uncertainty is revealed. Consequently, it solves the following min-max problem:
\vspace{-0.3cm}
\begin{subequations}\label{SRO}
\begin{align}[left={(SRO)}\empheqlbrace]
\nonumber \underset{x\in \Phi, y}{\min} \ &    \underset{u\in \mathbf{U}}{\max} \enspace    C_{I}^{\top}x + C_{O}^{\top}y\\
  s.t.     \ \enspace            & h(x,y,u) = 0  \quad  \label{Tri-Level-problem} \\ 
               & g(x,y,u) \leq 0  \quad  \\
\nonumber               & y \in \mathbb{R}
\end{align}
\end{subequations}
\noindent where $\mathbf{U}$  represents the uncertainty set of all realizations of the uncertain parameters $u$.  For the OSP-DER, the investment variables $x$ and operation variables  $y$ are simultaneously decided considering the worst-case value of the power generation and electricity demand. This approach can lead to overly conservative solutions and is therefore used in this study as an upper bound in the experiments.

\subsection{Adaptive Robust Optimization (ARO) formulation}\label{ch:chap4}

ARO follows a two-stage structure that considers second-stage variables that can be decided in response to the worst-case realization of uncertainty. This adaptability is achieved through the following  three-level optimization model:
\vspace{-0.3cm}
\begin{subequations}\label{ARO}
\begin{align}[left={(ARO)}\empheqlbrace]
\nonumber \underset{x\in \Phi}{\min} \;  C_{I}^{\top}x +     \underset{u\in \mathbf{U}}{\max} \; \underset{y} {\min}  & \enspace    C_{O}^{\top}y\\
  s.t.      \enspace            & h(x,y,u) = 0    \label{aroH}\\ 
               & g(x,y,u) \leq 0   \label{aroG}    \\
\nonumber               & y \in \mathbb{R}
\end{align}
\end{subequations}

For the OSP-DER, the first level determines the optimal placement and the capacity of DERs, the second level models the worst-case realization of the power generation and electricity demand,  and the third level minimizes the impact of the operational costs of the investment decisions. Moreover, since power systems operate within market structures where operational decisions significantly influence overall performance, ARO provides a more flexible and realistic framework for handling uncertainty and improving planning decisions.



\subsection{Adaptive Robust Stochastic Optimization (ARSO) formulation}

ARO approach assumes that all uncertainties exhibit the same nature, which may not always be true. For example, electricity demand uncertainty is predominantly an LT factor, influenced by structural changes in consumption patterns, policy shifts, and emerging technologies such as electric vehicles. Conversely, PV generation uncertainty arises from ST variability, primarily driven by weather conditions, which fluctuate hourly or daily. Therefore, treating both uncertainties identically may lead to overly conservative solutions. In this context, the ARSO approach combines the strengths of robust optimization and stochastic optimization to provide a more accurate representation of the distinct nature of uncertainties. The new hybrid three-level optimization model can be formulated  as follows:
\vspace{-0.5cm}
\begin{subequations}
\begin{align}[left={(ARSO)}\empheqlbrace]
 \nonumber \underset{x\in \Phi}{\min} \enspace  C_{I}^{\top}x +   \underset{u\in \mathbf{U}}{\max} \enspace \mathbb{E}_{\xi} [Q(x,&u,\xi)] \\
  \nonumber  \text{where} \quad  Q(x,u,\xi) = \min_{y} \enspace &  C_{O}^{\top}y \\
   s.t. \enspace  &h(x,y,u,\xi) = 0 \\
    &g(x,y,u,\xi) \leq 0 \\
   \nonumber  & y \in \mathbb{R} 
\end{align}
\end{subequations}

For the OSP-DER, the first level is the same as ARO of investment decisions, in contrast, the second-level problem only accounts for the worst-case realization of LT electricity demand within the predefined uncertainty set, and in the third-level the PV generation is represented as a random variable $\xi$, capturing its stochastic nature through a set of representative scenarios. This approach ensures resilience against LT electricity demand uncertainty and introduces adaptability in ST operational decisions related to PV generation uncertainty.



\subsection{Budgeted uncertainty sets}

The choice of the uncertainty sets type significantly influences the solution's conservatism and problem tractability. The simplest type is to use box uncertainty sets, i.e., $u\in \mathbf{U}  = [ \bar{u} - \hat{u}, \bar{u} + \hat{u} ]$, where $\bar{u}$ and  $\hat{u}$ are the expected value and deviation of each parameter, respectively. However, this can  lead to overly conservative solutions because it assumes that all uncertain parameters can simultaneously reach their worst-case values, which is often unrealistic. To mitigate this conservatism, the budgeted uncertainty sets are introduced to limit the number of simultaneous parameters that can deviate from their worst-case realizations given a budget of uncertainty $\beta$. By tuning this budget, the trade-off between robustness and cost-effectiveness can be adjusted, leading to less pessimistic solutions while still ensuring operational feasibility. Following~\cite{attarha2018adaptive}, we initially consider the following budgeted uncertainty set $\mathbf{U}$ for our uncertain parameters of power generation and electricity demand:
 

\begin{flalign}
\begin{split}
    \mathbf{U} = \bigg\{ (\widetilde{PV}_{t},\widetilde{PL}_{i,t})\in \mathbb{R^+} \enspace \bigg|   \enspace &|\widetilde{PV}_{t} - \overline{PV}_{t}| \leq \widehat{PV}_{t}, \enspace |\widetilde{PL}_{i,t} - \overline{PL}_{i,t}| \leq \widehat{PL}_{i,t},  \\ 
    &\displaystyle\sum_{i\in {\mathcal {B}}} \displaystyle\sum_{t\in {\mathcal {T}}} \bigg|\frac{\widetilde{PV}_{t} - \overline{PV}_{t}}{\widehat{PV}_{t}}\bigg| + \bigg|\frac{\widetilde{PL}_{i,t} - \overline{PL}_{i,t}}{\widehat{PL}_{i,t}}\bigg| \leq \beta \enspace;  \enspace \forall i \in {\mathcal {B}}, t \in {\mathcal {T}}\bigg\} \label{Gamma_Eq}
\end{split}
\end{flalign}

In the Set \eqref{Gamma_Eq}, $\widetilde{PV}_{t}$ and $\widetilde{PL}_{i,t}$ are individually defined in the intervals \mbox{$[\overline{PV}_{t}- \widehat{PV}_{t}, \overline{PV} _{t}+ \widehat{PV}_{t}]$} and \mbox{$[\overline{PL}_{i,t}- \widehat{PL}_{i,t}, \overline{PL} _{i,t}+ \widehat{PL}_{i,t}]$}, respectively, and at the same time the number of uncertain parameters that can simultaneously reach their worst-case deviations is limited by $\beta$. Thus, if $\beta$ is set to zero, both uncertainty parameters take their forecasted values (deterministic problem), whereas increasing $\beta$ raises robustness by allowing more uncertainty parameters to deviate from their expected values. If $\beta$ takes its maximum value, the model yields results equivalent to the original box uncertainty sets.

The BoU provides the decision-maker with a mechanism to regulate the model's robustness by increasing or reducing $\beta$. However, a single parameter $\beta$ in \eqref{Gamma_Eq}, may not adequately capture the distinct nature of different uncertainty sources. In particular, PV generation and electricity demand exhibit different variability patterns and levels of predictability, requiring independent control over their respective robustness levels. To address this, we propose decomposing the BoU into two separate parameters, $\beta_{pl}$ and $\beta_{pv}$, allowing for differentiated regulation of conservatism as follows:

\begin{subequations} \label{PV_PL_Polyhedral_Set}
\begin{flalign}
    &\displaystyle\sum_{t\in {\mathcal {T}}}\bigg|\frac{\widetilde{PL}_{i,t} - \overline{PL}_{i,t}}{\widehat{PL}_{i,t}}\bigg| \leq \beta_{pl} \enspace \enspace \forall i \in \mathcal {B} \label{PL_Polyhedral_Set}\\
    &\displaystyle\sum_{t\in {\mathcal {T}}} \bigg|\frac{\widetilde{PV}_{t} - \overline{PV}_{t}}{\widehat{PV}_{t}}\bigg| \leq \beta_{pv}
\end{flalign}
\end{subequations}

Note that Equation \eqref{PL_Polyhedral_Set} establishes a BoU for each bus, indicating that every bus with demand can deviate from its expected value up to $\beta$ times. However, caution is needed when applying this formulation, as including a large number of consumption points can significantly increase the computational complexity of the model. In such cases, defining a uniform BoU across all consumption points may be more efficient than differentiating by individual demand locations.

\section{Solving the ARO and ARSO OSP-DER formulation with BoU}\label{sec:algo}

The three-level ARO and ARSO formulations are not tractable for commercial optimization solvers due to their hierarchical structure. Typically, the second and third levels are merged into a single maximization problem~\cite{jeyakumar2010strong}, which can be efficiently handled using Benders decomposition, as commonly applied in the literature~\cite{bertsimas2012adaptive}. This section introduces the general framework of Benders decomposition algorithms and discusses strategies for effectively managing the budgeted uncertainty sets to solve the ARO and ARSO formulations.

\subsection{Benders decomposition algorithms }

The ARO and ARSO models introduced in the previous section can be solved using the Benders cutting plane algorithm~\citep{rahmaniani2020benders, ZENG2013457}. In a classical Benders decomposition, the original problem is split into a master problem (MP) and one subproblem (SP), which depends on the MP solution.  Thus, the original problem can be reformulated through Benders cuts that are built from the solutions and the objective function of the dual subproblem (DSP). Consequently, the original problem can be iteratively solved by solving the MP, which proposes a solution and gives a lower bound, as it is a relaxed problem, and then the DSP gives the corresponding upper bound. This type of cutting plane algorithm reduces the computational complexity to solve large-scale two-step problems by guaranteeing convergence in a finite number of iterations~\citep{thiele2009robust, zugno2015robust}.

Several strategies have been proposed in the literature to enhance the implementation of Benders decomposition algorithms~\citep{rahmaniani2020benders}. In this study, we adopt a specific variant where the MP initially uses the expected values for variables modeled under the robust optimization approach~\citep{zugno2015robust}. This modification provides a more informed starting point, improving solution quality and reducing the number of iterations required for convergence.

For the ARO formulation, the MP determines the capacity to be installed based on optimizing the first-stage variables $x$, while the subproblems, given this capacity as a fixed parameter $\hat{x}$, minimize the operational costs under the corresponding worst-case scenario. Note that the SP contains binary variables $w_{i,t}$ associated with battery charging and discharging, preventing the retrieval of dual variables needed to apply strong duality to the second- and third-level problems. To address this, we include the operational constraints related to these decisions in the MP, avoiding the need to linearize battery constraints, which would otherwise lead to unrealistic battery operation patterns. Consequently, the  MP and DSP are defined  as follows:
\vspace{-0.4cm}
\begin{align}[left={({MP})}\empheqlbrace]
 \nonumber   \underset{x\in \Phi}{\min} \enspace & C_{I}^{\top}x + \eta \label{MP_OB}\\
 s.t \enspace & h(x,y,{u}^{*}) = 0 \\
 \nonumber       &g(x,y,{u}^{*}) \leq 0\\
\nonumber    &\eta \geq [\mathbf{h}(x,{u}^{*})]^{\top} {\lambda}^{*} + [\mathbf{g}(x,{u}^{*})]^{\top} {\mu}^{*}
\end{align}
\begin{align}[left={({DSP(\hat{x})})}\empheqlbrace]
\nonumber\underset{{u},\lambda, \mu}{\max} \enspace & [\mathbf{h}(\hat{x},{u})]^{\top} \lambda + [\mathbf{g}(\hat{x},{u})]^{\top} \mu  \label{SP_OB}\\
 s.t. \enspace & [\mathbf{H}(\hat{x},{u})]^{\top} \lambda + [\mathbf{G}(\hat{x},{u})]^{\top} \mu = C_{O}\\
\nonumber&\lambda \in \mathbb{R}, \enspace \mu \in \mathbb{R}, \enspace  u \in \mathbf{U}
\end{align}
In the MP, ${u}^{*}$, $ {\lambda}^{*}$, and ${\mu}^{*}$ represent the worst-case realization of the uncertain parameters and solution obtained in the corresponding cut from the DSP at each iteration, by the auxiliary variable $\eta$. Initially, in the first iteration, ${u}^{*}$ is set to the expected value of the uncertainty parameter $\bar{u}$. In the DSP, the worst-case realization of ${u}$, is identified by maximizing the dual variables $\lambda$ and $\mu$, which enforce the operational constraints of the primal inner problem while ensuring that the level of adversity remains within the predefined BoU  $\beta$ previously defined in $\mathbf{U}$. The resulting worst-case realization ${u}^{*}$ is then added to the current MP for the next iteration, iteratively refining the solution through the Benders cuts. The Benders decomposition algorithm for the ARO formulation is summarized in Algorithm~\ref{Algorithm_label_ARO}.

\vspace{0.3cm}
\begin{algorithm}[H]
    \caption{Benders decomposition algorithm for ARO formulation}
    \label{Algorithm_label_ARO}
    \SetKwInOut{Init}{Initialize}
    \Init{
       $LB \leftarrow -\infty, \; UB \leftarrow + \infty, \; r \leftarrow 0$\\
    }
    \While{$UB > LB$}{ 
        $(\widehat{x},\, \widehat{y})  \gets $ Solve MP \\
         $LB = C_{I}^{\top}\widehat{x} +  \eta$\\
            $({\lambda}^{*}, \, {\mu}^{*}, \,  \theta ) \gets$    Solve $DSP(\widehat{x})$  \\
            Add the Benders cut $\eta \geq [\mathbf{h}(x,{u}^{*})]^{\top} {\lambda}^{*} + [\mathbf{g}(x,{u}^{*})]^{\top} {\mu}^{*}$ to the MP \\
    	   $UB \gets \min \big(UB,\,  LB -  (\eta + \theta) \big)$\\
    }
        \Return {$(\widehat{x},\,  \widehat{y},\,   UB)$}

\end{algorithm}
\vspace{0.3cm}

In Steps 2 and 3, we compute the MP solution of the first-stage decisions $\widehat{x}$, the corresponding operational decision $\widehat{y}$, and the objective value as the lower bound $LB$. In Step 4, we obtain the dual variables $({\lambda}^{*}, {\mu}^{*})$, and its objective value $\theta$. Then, in Step 5, the corresponding Benders cut is added to the MP, and the upper bound $UB$ is updated in Step 6. The process is repeated if there is still an optimality gap.

For the ARSO formulation, it is possible to generate a Benders cut for each scenario at each iteration, specifically for PV generation decision variables. This multi-cut strategy enhances the algorithm's convergence by incorporating multiple cuts into the MP in each iteration, improving computational efficiency~\cite{GARCIAMUNOZ2023120226}. Consequently, the corresponding MP and DSP follow a similar structure but include the subindex of scenarios $s\in S$ in the corresponding Benders cuts and then use $|S|$ auxiliary variables $\eta_s$. The Benders decomposition algorithm for the ARSO formulation is summarized in Algorithm~\ref{Algorithm_label_ARSO}.

\vspace{0.3cm}
\begin{algorithm}[H]
    \caption{Benders decomposition algorithm for ARSO formulation}
    \label{Algorithm_label_ARSO}
    \SetKwInOut{Init}{Initialize}
    \Init{
       $LB \leftarrow -\infty, \; UB \leftarrow + \infty, \; r \leftarrow 0$\\
    }
    \While{$UB > LB$}{ 
        $(\widehat{x},\, \widehat{y})  \gets $ Solve MP \\
          $LB = C_{I}^{\top}x + \displaystyle\sum_{s\in {\mathcal {S}}} \eta_s$\\
          \For{$s \in S$}{
            $({\lambda}^{*}, \,  {\mu}^{*}, \,  \theta_s ) \gets$    Solve $DSP(\widehat{x})$  \\
            Add the Benders cut $\eta_s \geq [\mathbf{h}(x,{u}^{*},s)]^{\top} {\lambda}^{*} + [\mathbf{g}(x,{u}^{*},s)]^{\top} {\mu}^{*}$ to the MP \\
    	   $UB \gets \min \big(UB,\,  LB -  \displaystyle\sum_{s\in {\mathcal {S}}} (\eta_s + \theta_s) \big)$\\
          }
  
    }
            \Return {$(\widehat{x},\,  \widehat{y},\,   UB)$}

\end{algorithm}
\vspace{0.3cm}

Similarly to Algorithm 1, in Steps 2 and 3, we compute $\widehat{x}$, $\widehat{y}$, and $LB$, now with the auxiliary variables $\eta_s$. For each scenario $s\in S$, the corresponding Benders cuts and $UB$ are computed in Steps 4 to 7. This procedure is repeated until an optimal solution is found.

\subsection{Dealing with the BoU within the dual sub-problems}

The budgeted uncertainty sets introduced in \eqref{Gamma_Eq} define a polyhedral feasible region for $\widetilde{PV}_{t}$ and $\widetilde{PL}_{i,t}$, where the optimal solution lies at one of its extreme points~\cite{attarha2018adaptive}. Following~\cite{jiang2011robust}, for each $i \in {\mathcal {B}}$ and $t \in {\mathcal {T}}$, new binary variables $U^{+}_{t}, U^{-}_{t}, V^{+}_{i,t}, V^{-}_{i,t}$ are incorporated into the model to effectively identify these extreme points via the following constraints:
\begin{subequations}\label{BOU_Constraints}
\begin{flalign}
    &\widetilde{PV}_{t} =  \overline{PV}_{t} + \widehat{PV}_{t}U^{+}_{t} - \widehat{PV}_{t}U^{-}_{t}, \label{BOU_Constraints_1}\\
    &\widetilde{PL}_{i,t} =  \overline{PL}_{i,t} + \widehat{PL}_{i,t}V^{+}_{i,t} - \widehat{PL}_{i,t}V^{-}_{i,t} \label{BOU_Constraints_2}\\
    &\displaystyle\sum_{t\in {\mathcal {T}}} (U^{+}_{t} + U^{-}_{t}) \leq \beta_{pv}\label{BOU_Constraints_3} \\
    &\displaystyle\sum_{t\in {\mathcal {T}}} (V^{+}_{i,t} + V^{-}_{i,t}) \leq \beta_{pl} \label{BOU_Constraints_4}\\
    &U^{+}_{t} + U^{-}_{t} \leq 1 \label{BOU_Constraints_5}\\
    &V^{+}_{i,t} + V^{-}_{i,t} \leq 1 \label{BOU_Constraints_6}
\end{flalign}
\end{subequations}

\vspace{0.4cm}

Constraints \eqref{BOU_Constraints_1} and \eqref{BOU_Constraints_2} model these deviations using the binary variables, ensuring that each uncertain parameter can only shift in one direction within the uncertainty budget constraints, as enforced by \eqref{BOU_Constraints_5} and \eqref{BOU_Constraints_6}. The budgets of uncertainty in Constraints \eqref{BOU_Constraints_3} and \eqref{BOU_Constraints_4}  determine  the number of extreme deviations. For instance, if $\beta_{pv} = 5$ in a 24-hour horizon with an hourly time step, Constraints \eqref{BOU_Constraints_3} ensures that $\widetilde{PV}_{t}$ can deviate from the expected value $\overline{PV}_{t}$ in at most five out of the 24 forecasted hours. These deviations occur through the activation of either $U^{+}_{t}$ or $U^{-}_{t}$, depending on the worst-case scenario, but never simultaneously given Constraint \eqref{BOU_Constraints_5}. Similarly, the same logic applies for $\widetilde{PL}_{i,t}$ enforced for each $i \in {\mathcal {B}}$. For example, considering $\beta_{pl} = 5$, each bus can undergo up to five deviations from its forecast demand value.

Note that nonlinearity is introduced in the corresponding objective functions by incorporating Constraints \eqref{BOU_Constraints} in the corresponding DSP \eqref{SP_OB}. To explicitly show this, the following is the extended objective function of the DSP in the ARO formulation: 

\vspace{-0.6cm}

\begin{multline}\label{ARO_Extensive_Formulation}
\underset{\lambda, \mu }{\max} \displaystyle\sum_{i\in {\mathcal {B}}} \displaystyle\sum_{t\in {\mathcal {T}}} \bigg\{   -\widetilde{PL}_{i,t}\mathbf{a}_{i,t} - QL_{i,t}\mathbf{b}_{i,t} + 
PG^{max}_{i}\mathbf{f}_{i,t} + 
\mathbf{g}^{min}_{i,t}QG^{min}_{i} +  \mathbf{g}^{max}_{i,t}QG^{max}_{i}  
+ \mathbf{h}^{min}_{i,t}V^{min}_{i} \\ +  \mathbf{h}^{max}_{i,t}V^{max}_{i} + 
\widetilde{PV}_{t} \widehat{\gamma}^{pv}_{i}\mathbf{i}_{i,t} + \mathbf{k}^{min}_{i,t}SOC^{min}\widehat{\gamma}^{bt}_{i} 
+  \mathbf{k}^{max}_{i,t}SOC^{max} \widehat{\gamma}^{bt}_{i}  + 
\mathbf{l}_{i,t}PB\widehat{w}_{i,t}  \\ + \mathbf{m}_{i,t}(PB(1-\widehat{w}_{i,t})) +  \mathbf{n}^{min}_{i,t}\Theta^{min}_{i} + 
\mathbf{n}^{max}_{i,t}\Theta^{max}_{i}     \bigg\}  
- \displaystyle\sum_{(i,j)\in {\mathcal {L}}} \displaystyle\sum_{r\in {\mathcal {R}}} S^{max}_{i,j}C_{r}\mathbf{e}^{r}_{i,j,t} + \displaystyle\sum_{i\in {\mathcal {B}}}  SOC^{init} \widehat{\gamma}^{bt}_{i} \mathbf{j}_{i,1} 
\end{multline}

\noindent where the dual variables associated with the equality and inequality constraints are $\lambda = \{\mathbf{a}_{i,t},\mathbf{b}_{i,t}, \mathbf{j}_{i,1}   \}$ and $\mu = \{\mathbf{f}_{i,t},\mathbf{g}^{min}_{i,t},\mathbf{g}^{max}_{i,t},\mathbf{h}^{min}_{i,t},\mathbf{h}^{max}_{i,t},\mathbf{i}_{i,t}, \mathbf{k}^{min}_{i,t},\mathbf{k}^{max}_{i,t},\mathbf{l}_{i,t},\mathbf{m}_{i,t}, \mathbf{n}^{min}_{i,t}, \mathbf{n}^{max}_{i,t},\mathbf{e}^{r}_{i,j,t}    \}$, respectively. The respective dual variables corresponding to each constraint of the subproblem are specified in~\ref{appendix:spARO}, and the extended dual subproblem can be found in~\ref{appendix:dspARO}. The nonlinearity in \eqref{ARO_Extensive_Formulation} arise from the product of dual variables $\mathbf{a}_{i,t}$ and $\mathbf{i}_{i,t}$ with the uncertain parameters $\widetilde{PL}_{i,t}$ and $\widetilde{PV}_{t}$, respectively. These bilinearities make the problem non-convex, increasing computational complexity. To address this issue, a Big-M method \cite{attarha2018adaptive} is employed by introducing two sets of auxiliary continuous variables for each dual variable, $(\mathbf{a}^{+}, \mathbf{a}^{-})$ and $(\mathbf{i}^{+}, \mathbf{i}^{-})$ for each $ i \in {\mathcal {B}}$ and $ t \in {\mathcal {T}}$, as follows:

\begin{subequations}\label{BigM_Eqs}
\small
\begin{flalign}
    &\widetilde{PL}_{i,t}\mathbf{a}_{i,t}  = \overline{PL}_{i,t}\mathbf{a}_{i,t} + \widehat{PL}_{i,t}(\mathbf{a}^{+}_{i,t} - \mathbf{a}^{-}_{i,t})\label{PL_bilinear_termn2}\\
    &\widetilde{PV_{t}}\widehat{\gamma}^{pv}_{i}\mathbf{i}_{i,t}  = \overline{PV}_{t}\widehat{\gamma}^{pv}_{i}\mathbf{i}_{i,t} + \widehat{PV}_{t}(\mathbf{i}^{+}_{i,t} - \mathbf{i}^{-}_{i,t})\widehat{\gamma}^{pv}_{i}\label{PV_bilinear_termn2}\\
    &-V^{+}_{i,t} M \leq \mathbf{a}^{+}_{i,t} \leq V^{+}_{i,t} M \\
    &-V^{-}_{i,t} M \leq \mathbf{a}^{-}_{i,t} \leq V^{-}_{i,t} M \\
    &-U^{+}_{t} M \leq \mathbf{i}^{+}_{i,t} \leq U^{+}_{t} M \\
    &-U^{-}_{t} M \leq \mathbf{i}^{-}_{i,t} \leq U^{-}_{t} M \\
    &\mathbf{a}_{i,t} - M(1-V^{+}_{i,t}) \leq \mathbf{a}^{+}_{i,t} \leq \mathbf{a}_{i,t} + M(1-V^{+}_{i,t}) \\
    &\mathbf{a}_{i,t} - M(1-V^{-}_{i,t}) \leq \mathbf{a}^{-}_{i,t} \leq \mathbf{a}_{i,t} + M(1-V^{-}_{i,t}) \\
    &\mathbf{i}_{i,t} - M(1-U^{+}_{t}) \leq \mathbf{i}^{+}_{i,t} \leq \mathbf{i}_{i,t} + M(1-U^{+}_{t}) \\
    &\mathbf{i}_{i,t} - M(1-U^{-}_{t}) \leq \mathbf{i}^{-}_{i,t} \leq \mathbf{i}_{i,t} + M(1-U^{-}_{t}) 
\end{flalign}
\end{subequations}

Finally, introducing the auxiliary variables in Constraints \eqref{BigM_Eqs} ensures a linear objective function of the dual subproblem \eqref{ARO_Extensive_Formulation}. Note that $M$ must be large enough to cover the feasible range of uncertain parameters without excessive conservatism.  It can be computed as: 

\begin{equation}
    M = \max \left\{ \max_{i \in \mathcal{B}, t \in \mathcal{T}} (\overline{PL}_{i,t} + \widehat{PL}_{i,t}), \, \max_{t \in \mathcal{T}} (\overline{PV}_{t} + \widehat{PV}_{t}) \right\} 
\end{equation}

The ARSO DSP follows a similar approach to the DSP of the ARO formulation,  with the difference in handling PV generation uncertainty through discrete scenarios. The extended  DSP objective function of ARSO is as follows:

\vspace{-0.4cm}
\begin{multline}\label{ARSO_Extensive_Formulation}
\underset{\lambda, \mu }{\max} \displaystyle\sum_{i\in {\mathcal {B}}} \displaystyle\sum_{t\in {\mathcal {T}}} \displaystyle\sum_{s\in {\mathcal {S}}}\bigg\{   -\widetilde{PL}_{i,t}\mathbf{a}_{i,t,s} - QL_{i,t,s}\mathbf{b}_{i,t,s} + 
PG^{max}_{i}\mathbf{f}_{i,t,s} + 
\mathbf{g}^{min}_{i,t,s}QG^{min}_{i} +  \mathbf{g}^{max}_{i,t,s}QG^{max}_{i} \\ 
+ \mathbf{h}^{min}_{i,t,s}V^{min}_{i} +  \mathbf{h}^{max}_{i,t,s}V^{max}_{i} + 
PV_{t,s} \widehat{\gamma}^{pv}_{i}\mathbf{i}_{i,t,s} + \mathbf{k}^{min}_{i,t,s}SOC^{min}\widehat{\gamma}^{bt}_{i} 
+  \mathbf{k}^{max}_{i,t,s}SOC^{max} \widehat{\gamma}^{bt}_{i} \\ + 
\mathbf{l}_{i,t,s}PB\widehat{w}_{i,t,s} + \mathbf{m}_{i,t,s}(PB(1-\widehat{w}_{i,t,s})) +  \mathbf{n}^{min}_{i,t,s}\Theta^{min}_{i} + 
\mathbf{n}^{max}_{i,t,s}\Theta^{max}_{i}     \bigg\} \\ 
- \displaystyle\sum_{(i,j)\in {\mathcal {L}}} \displaystyle\sum_{r\in {\mathcal {R}}} \displaystyle\sum_{s\in {\mathcal {S}}} S^{max}_{i,j}C_{r}\mathbf{e}^{r}_{i,j,t,s} + \displaystyle\sum_{i\in {\mathcal {B}}} \displaystyle\sum_{s\in {\mathcal {S}}} SOC^{init} \widehat{\gamma}^{bt}_{i} \mathbf{j}_{i,1,s} 
\end{multline}

The  extended ARSO dual subproblem  can be found in~\ref{appendix:dspARSO}. Note that in \eqref{ARSO_Extensive_Formulation}, only $\widetilde{PL}_{i,t}$ is modeled using BoU, while $PV_{t,s}$ is represented through discrete scenarios indexed by $s$. Therefore, the constraints associated with PV generation in \eqref{BigM_Eqs} must not be included in the ARSO DSP. The bilinearities related to $\widetilde{PL}_{i,t}$ are handled similarly to \eqref{BigM_Eqs}, except that the dual variables $\mathbf{a}_{i,t,s}$ and the auxiliary variables $\mathbf{a}^{+}_{i,t,s}$ and $\mathbf{a}^{-}_{i,t,s}$ must include the scenario index $s$.

\newpage
\section{Computational results}\label{sec:results}

This section presents the case study used to evaluate how the proposed formulations handle LT and ST uncertainties and their impact on the installed capacity of DERs in terms of the conservatism and autonomy levels of the corresponding solutions. We used Pyomo 6.7.3 and CPLEX 12.7.0 on a laptop equipped with an Intel Core i9 processor, 24 cores, and 64 GB of memory.

\subsection{Case Study}
A modified version of the well-known IEEE 33-bus system~\cite{dolatabadi2020enhanced} was used as the test system to evaluate the performance of the stochastic and robust formulations presented in the previous section. It consist of 33 buses and 32 branches, where 32 buses have power loads.  

Regarding the planning horizon, simulating a 20-year period with an hourly time step is computationally infeasible and inefficient~\cite{teichgraeber2022time}. Instead, the proposed methodology adopts representative periods to balance computational tractability with the ability to capture key operational patterns. This approach reduces the computational burden while preserving the accuracy of the capacity estimates. The operational costs calculated in the model are limited to these representative periods and do not account for the entire LT operational timeline. To ensure consistency, all capacity-related costs are appropriately adjusted to align with the representative periods used in the analysis. It is important to note that the scope of this article is restricted to estimating the capacity of DERs to be installed. A comprehensive assessment of the profitability of investment and operational costs over the full planning horizon is beyond the scope of this study. 

\begin{figure}[htpb]
\centering
\includegraphics[width=3.3in]{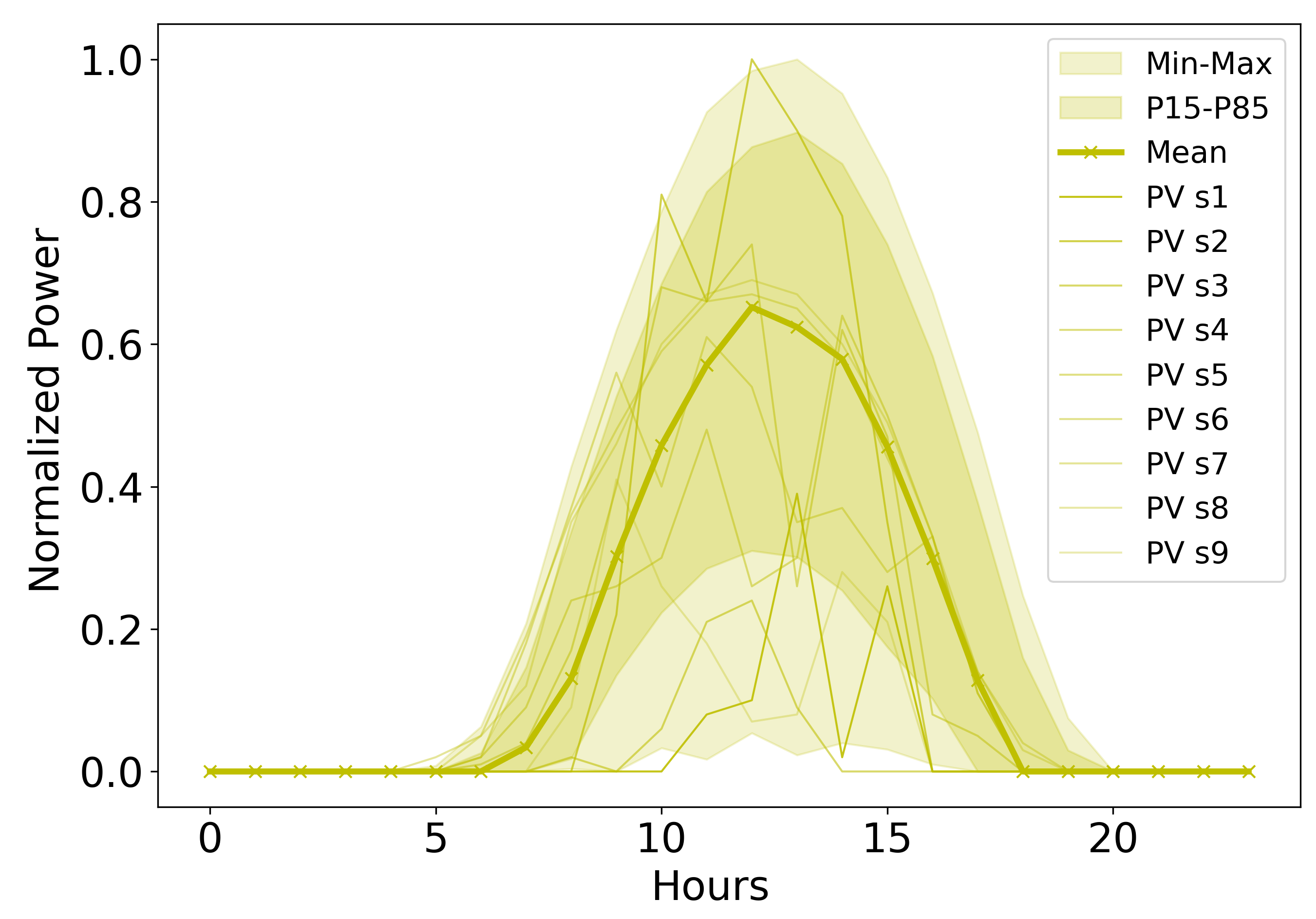}
\caption{PV Boundary set and scenarios}
\label{Plot_PV_Percentile}
\end{figure}

\begin{figure}
    \begin{subfigure}[b]{0.5\textwidth}
    \centering
    \includegraphics[width=3.3in]{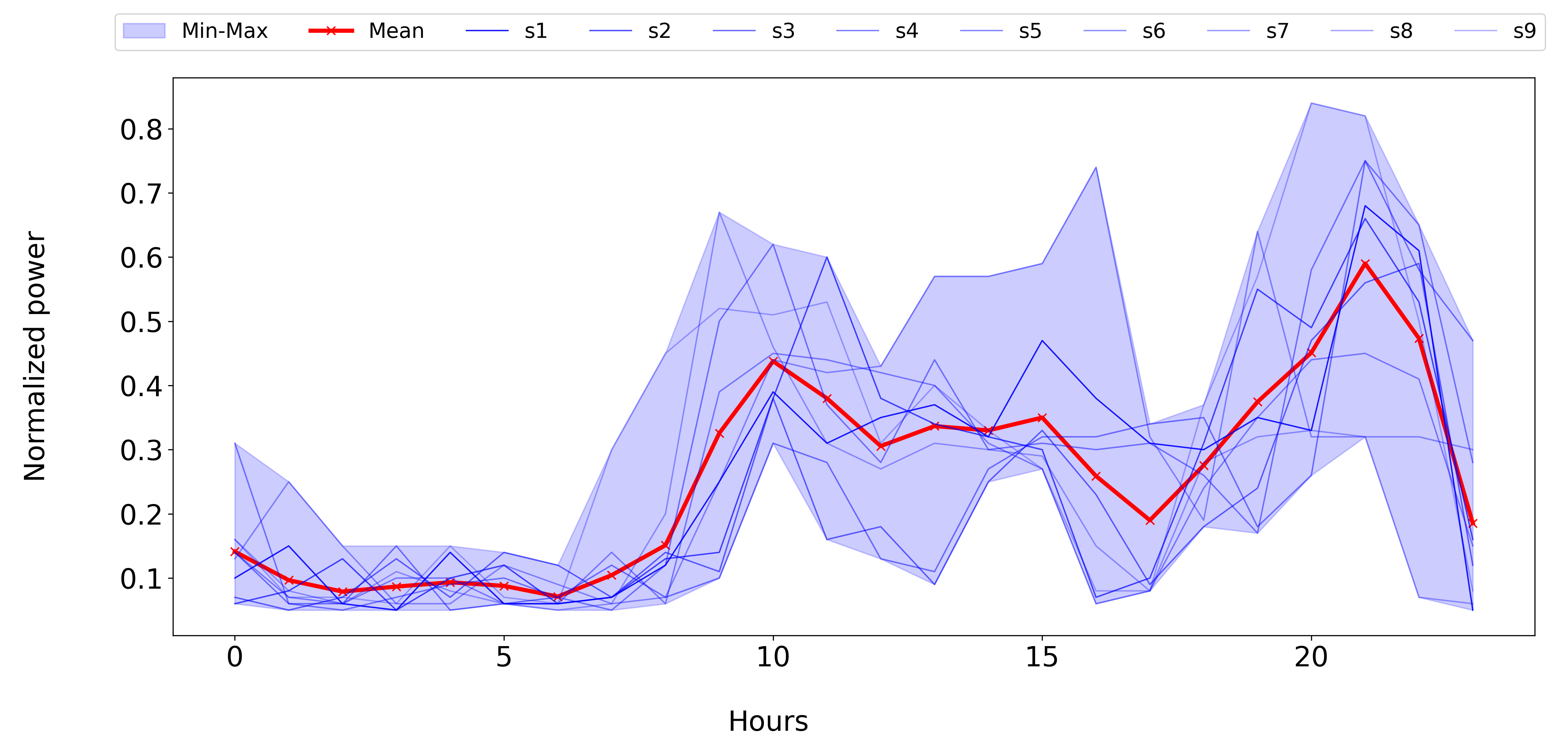}
    \caption{On bus 11}
    \label{Plot_Demand_Bus11}
    \end{subfigure}
    \begin{subfigure}[b]{0.5\textwidth}
        \centering
    \includegraphics[width=3.3in]{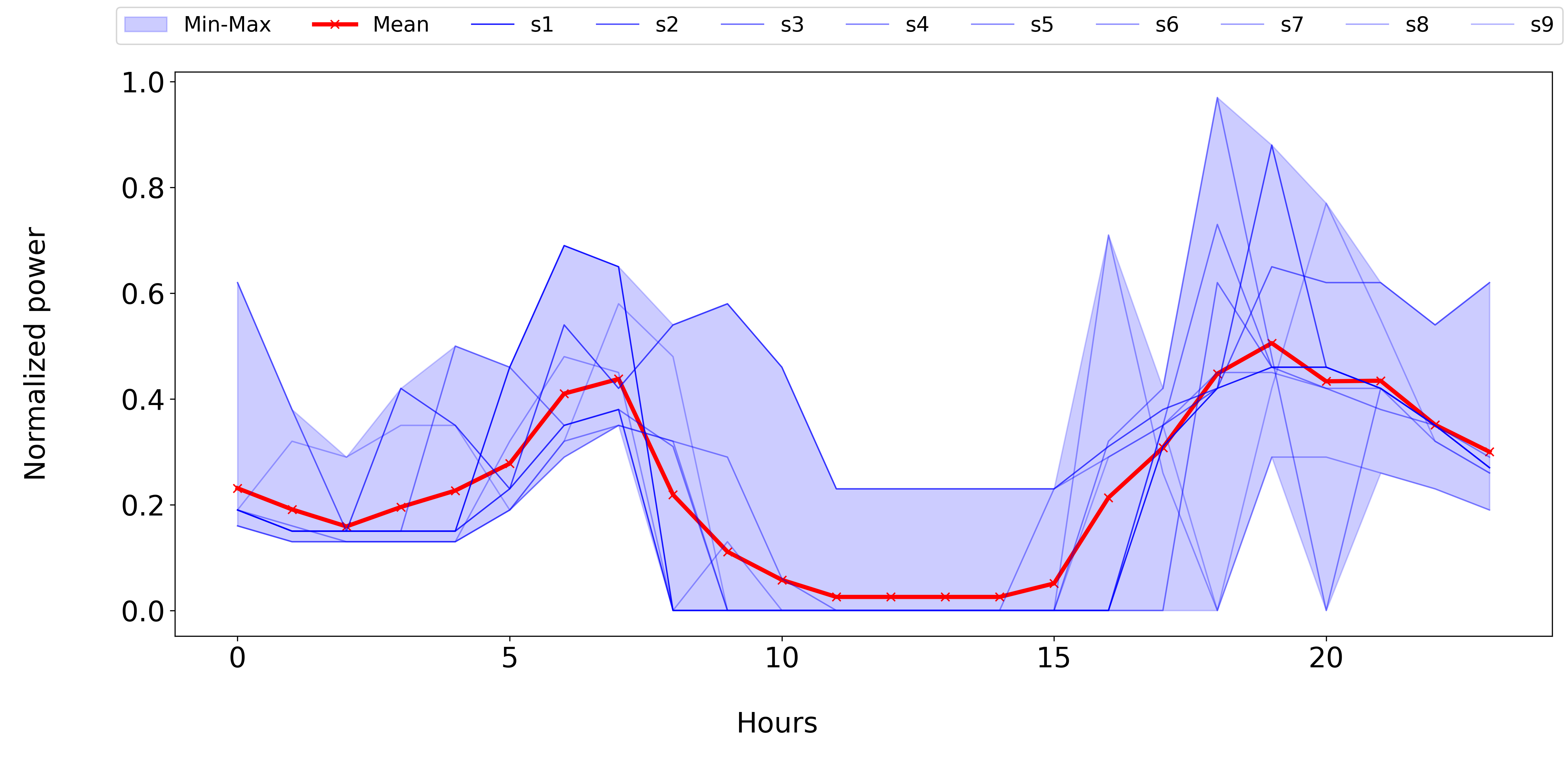}
    \caption{On bus 18}
    \label{Plot_Demand_Bus18}
    \end{subfigure}
    \caption{Boundary set and scenarios of electricity demand}
\end{figure}

Annual historical data were used to model the uncertainty in the hourly power injection of a PV system and electricity load consumption. Two distinct demand groups were considered: one representing traditional demand profiles without self-generation~\cite{Load_Data}, and another reflecting users equipped with PV systems, BESS, or both~\cite{GARCIAMUNOZ2024101538}. Figure~\ref{Plot_Demand_Bus11} illustrates a representative demand profile for the first group, while Figure~\ref{Plot_Demand_Bus18} depicts a profile for the second group. To model uncertainty as a discrete set of scenarios, as in the TSSO or ARSO approaches, nine representative scenarios were extracted from historical data using the backward reduction algorithm~\cite{heitsch2003scenario}, capturing the range of potential PV outputs. The nine electricity demand scenarios were derived directly from the irradiance scenarios selected by the backward algorithm, ensuring consistency between demand and generation conditions. For the ARO model with BoU, the expected value of the historical data was used as $\overline{PV}$, while the 85th and 15th percentiles were used as bounds for the maximum deviation $\widehat{PV}$ from the expected value. Figure~\ref{Plot_PV_Percentile} illustrates both the nine discrete scenarios and the BoU for the PV random parameter, highlighting the boundaries and variability of the data. This approach excluded extreme values below the 15th percentile to prevent overly pessimistic or optimistic scenarios, providing a realistic and probabilistic representation of uncertainty suitable for a 15–20-year planning horizon. Additionally, the polyhedral uncertainty set for load profiles was constructed based on the nine discrete scenarios selected, using their average as the expected value $\overline{PL}$ and the maximum and minimum consumption levels as bounds for the maximum deviation $\widehat{PL}$.

It became essential to establish a Perfect Information (PI) benchmark to properly evaluate how each formulation handles uncertainty. In this regard, PI refers to an ideal scenario where the future realization of uncertain parameters, such as electricity demand and PV generation, is fully known in advance. Although such a scenario is unrealistic in practice, it serves as a reference for assessing the trade-offs between robustness and cost-effectiveness in DER planning. A model that significantly deviates from the PI solution may either over-invest due to excessive conservatism or underestimate uncertainty, leading to unreliable system operation. To establish this benchmark, a 216-hour scenario was created by concatenating the nine representative scenarios selected via the backward reduction algorithm, providing a structured case to evaluate how closely different optimization models approximate the optimal solution under perfect foresight. Detailed information on network topology, parameters, operational and investment costs, and scenario data are available in the supplementary material.

\subsection{Reduced-scale testing}
The case study, based on 32 different consumption profiles, was initially conducted using a single load profile at a time to assess computational feasibility and identify preliminary patterns. In this exploratory analysis, the models were run separately for each load profile (A, B, and C) to determine the optimal sizing of PV systems alone, followed by a combined PV and BESS configuration. Profiles A and B represent traditional consumption without PV or BESS, while Profile C corresponds to a user with both PV and BESS installed. Additionally, the deterministic model was employed using the 216-hour operational scenario under PI to assess the impact of LT and ST uncertainty management across different optimization formulations. Table \ref{Table_Capacity_Profiles} presents the results, where the "PI" column denotes the optimal installed capacity under the 216-hour scenario, while the optimization formulation columns illustrate the deviation of each model from this optimal solution under perfect foresight.

\begin{table}[htpb]
\caption{Deviation in installed capacity across approaches relative to perfect information for different profiles and cases.}\label{Table_Capacity_Profiles}
\footnotesize
\resizebox{\textwidth}{!}{%
\begin{tabular}{cccccccccc}
\hline
Model               & Tech & Profile & PI & \textbf{$ARSO_{\beta=0}$} & \textbf{$ARO_{\beta=0}$} & \textbf{$ARSO_{\beta=10}$} & \textbf{$ARO_{\beta=10}$} & \textbf{$ARSO_{\beta=20}$} & \textbf{$ARO_{\beta=20}$} \\ \hline
\multirow{2}{*}{Only PV}     & PV            & A                & 0.053       & 0.001            & -0.003          & 0.005             & 0.014            & 0.048             & 0.084            \\
                             & PV            & B                & 0.050       & -0.050           & 0.000           & 0.005             & 0.020            & 0.027             & 0.087            \\ \hline
\multirow{4}{*}{PV and BESS} & PV            & A                & 0.053       & 0.004            & -0.002          & 0.012             & 0.034            & 0.050             & 0.089            \\
                             & BESS          & A                & 0.000       & 0.038            & 0.069           & 0.060             & 0.053            & 0.087             & 0.090            \\
                             & PV            & B                & 0.050       & 0.000            & 0.000           & 0.007             & 0.055            & 0.031             & 0.085            \\
                             & BESS          & B                & 0.000       & 0.000            & 0.000           & 0.040             & 0.058            & 0.051             & 0.036            \\ \hline
\end{tabular}
}
\end{table}

Profile C was not included in Table~\ref{Table_Capacity_Profiles} because none of the optimization models resulted in any PV or BESS installed capacity, resulting in zero capacity for both technologies. This is due to Profile C's intrinsic consumption pattern, which already incorporates self-generation capabilities. As a result, the cost-benefit trade-off did not justify additional PV and BESS installations. For the traditional consumption profiles A and B, the deviation from the PI column is lower for the ARSO model compared to the ARO formulation. This pattern holds for both consumption profiles and models that install only PV, as well as those that install PV plus BESS. Furthermore, as the BoU increases, the trend persists, with ARSO consistently exhibiting lower deviation than ARO.

\subsection{Conservatism in robust and stochastic models}

This section presents the full case study evaluation for the OSP-DER problem under different stochastic and robust formulations, including TSSO, SRO, ARO, and ARSO with varying BoU ($\beta$) levels in 0, 10, and 20, while considering the insights from the reduced-scale testing as a guiding reference. Table \ref{Table_Sizing_Allocation} summarizes the installed capacities of PV systems and BESS for each approach. Since robust optimization inherently seeks solutions that remain feasible under worst-case scenarios, more conservative models tend to allocate larger DER capacities to mitigate uncertainty. Thus, the level of conservatism in each formulation can be inferred from the total investment in DER capacity. The results in Table \ref{Table_Sizing_Allocation} highlight these differences, showing how each approach balances investment costs and uncertainty mitigation, while the specific allocation of technologies tends to remain consistent across methods. As expected, the SRO model incurs the highest investment cost in PV capacity due to less information being available at the time of decision-making. Similarly, the TSSO model installs PV but avoids BESS installation because its probabilistic approach does not prioritize extreme cases where BESS might be essential, and the cost-benefit analysis under expected scenarios makes batteries appear less economical compared to direct investments in PV capacity.

The ARSO and ARO approach effectively reduce the conservatism of the SRO across different $\beta$ values. Notably, ARSO further minimizes conservatism compared to ARO by differentiating between LT and ST uncertainties, whereas ARO, like SRO, models both uncertainties using the BoU approach but optimizes against the worst-case scenario to minimize operational costs. Specifically, as $\beta$ increases, the installed PV capacity in ARSO remains relatively stable. In contrast, the ARO model exhibits an increase in PV capacity with higher $\beta$, reflecting a more conservative response to uncertainty. However, the installed BESS capacity decreases as $\beta$ grows in both models. The ARO model consistently results in higher BESS capacity installations than ARSO, likely due to its stronger emphasis on robustness, cost-benefit trade-offs, and the utilization of surplus PV generation to meet demand.



\begin{table}[htpb]
\caption{OSP-DER for the different stochastic and robust optimization formulations}\label{Table_Sizing_Allocation}
\footnotesize
\resizebox{\textwidth}{!}{%
\begin{tabular}{ccccccccc}
\toprule
             \textbf{FORMULATION}      & \textbf{$TSSO$} & \textbf{$SRO$} & \textbf{$ARSO_{\beta=0}$} & \textbf{$ARSO_{\beta=10}$} & \textbf{$ARSO_{\beta=20}$} & \textbf{$ARO_{\beta=0}$} & \textbf{$ARO_{\beta=10}$} & \textbf{$ARO_{\beta=20}$} \\ \midrule
\textbf{Nº Bus PV} & 25              & 24             & 24                        & 25                         & 25                         & 24                       & 24                        & 25                        \\
\textbf{Nº Bus BESS} & -               & -              & 7                         & 7                          & -                          & 8                        & 24                        & 7                         \\  \midrule
\textbf{Cap. Pv [MW]}    & 1.68           & 3.79           & 1.59                      & 1.59                       & 1.66                       & 1.16                     & 1.86                      & 1.86                      \\
\textbf{Cap. BESS [MWh]}    & 0.00               & 0.00              & 0.18                      & 0.18                       & 0.00                       & 1.58                     & 0.97                      & 0.97                      \\ \midrule
\textbf{Inv. Cost [\euro]} & 293          & 658         & 311                    & 312                     & 287                      & 517                   & 517                    & 518                   \\
\textbf{Op. Cost [\euro]}  & 2569         & 2879        & 2527                   & 2881                    & 3154                     & 2292                  & 2753                   & 3002                   \\
\textbf{OF value [\euro]}      & 2862         & 3537        & 2839                   & 3193                    & 3441                     & 2809                  & 3271                   & 3519                  \\ \bottomrule
\end{tabular}
}
\end{table}

As shown in Table~\ref{Table_Sizing_Allocation}, the ARSO approach effectively reduces the conservatism of the SRO and ARO methods. However, the table does not quantify how its solution compares to the case where future conditions are perfectly known. To address this, we replicated the analysis from the previous section by running the deterministic model over the 216-hour scenarios, providing a benchmark to evaluate how far the ARSO and ARO solutions deviate from the installed capacity that would result under PI.

\subsection{Performance of ARO and ARSO models against perfect information}
\begin{figure}[htpb]
\centering
\includegraphics[width=4.0in]{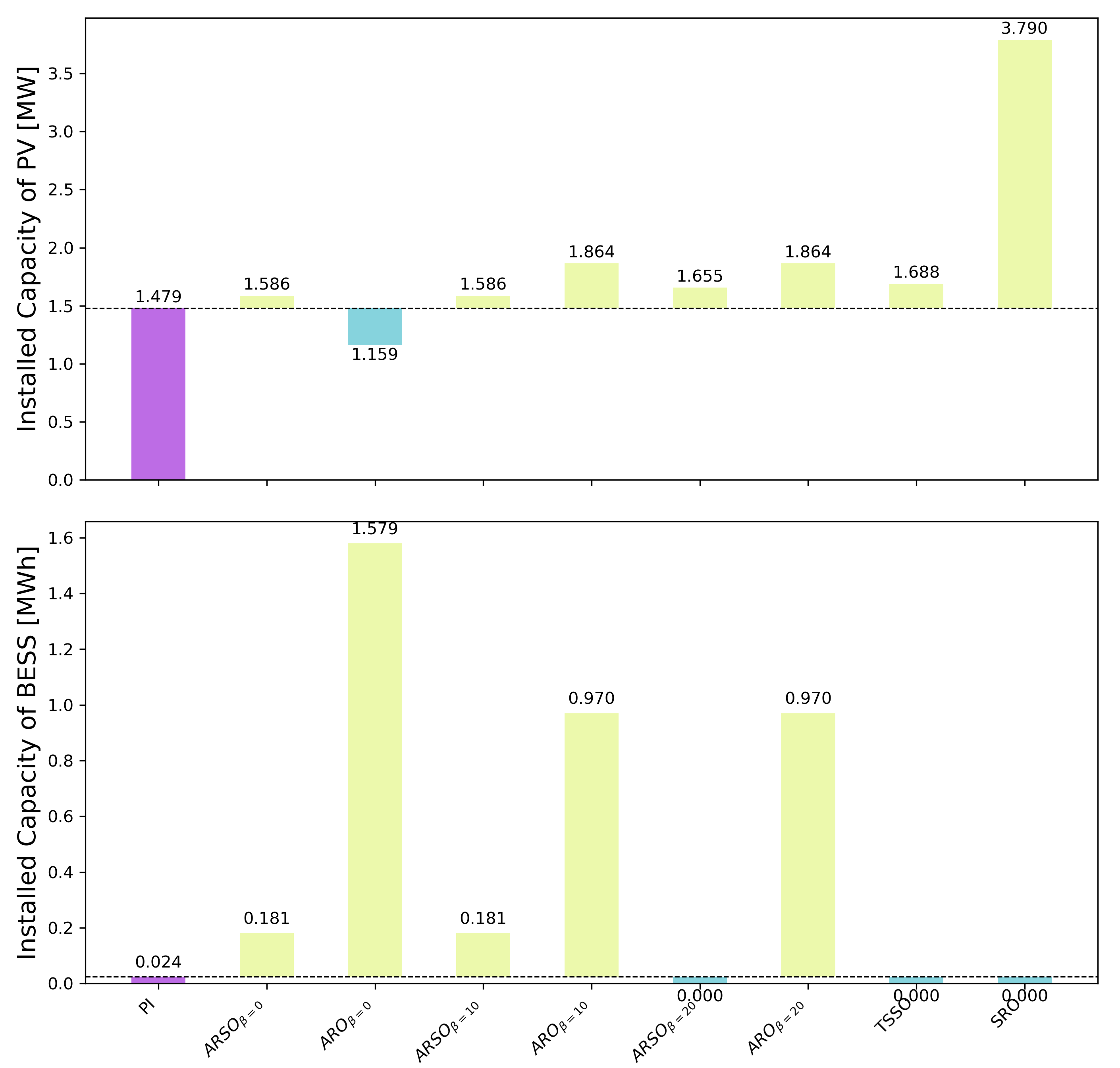}
\caption{Comparison of the installed capacities of PV and BESS across different approaches and the BoU, relative to PI.}
\label{Waterfall_Plot}
\end{figure}

Similar to the analysis conducted in the reduced-scale testing, this evaluation extends the comparison to the full case study, focusing on how ARO and ARSO models with different BoU values deviate from the optimal solution under PI. Figure~\ref{Waterfall_Plot} presents this comparison, where the purple bar represents the optimal installed capacity under a PI scenario, assuming full knowledge of future weather and demand conditions. The other bars illustrate the deviations from this reference value for PV systems (upper plot) and BESS (lower plot). Yellow bars indicate cases where the installed capacity exceeds the PI reference, while cyan bars represent cases where it falls below it. Notably, in the BESS plot, cyan bars appear when no batteries are installed. As observed, the ARSO approach achieves installed capacities closer to the PI reference compared to the ARO approach. This suggests that differentiating the treatment of uncertainties, as ARSO does, reduces conservatism while also resulting in solutions that align more closely with the ideal capacities under PI. In contrast, the ARO method, while robust, tends to exhibit greater deviations from the PI case, particularly by installing higher capacities, reflecting its inherently more conservative nature.

\begin{figure}[htpb]
\centering
\includegraphics[width=4.8in]{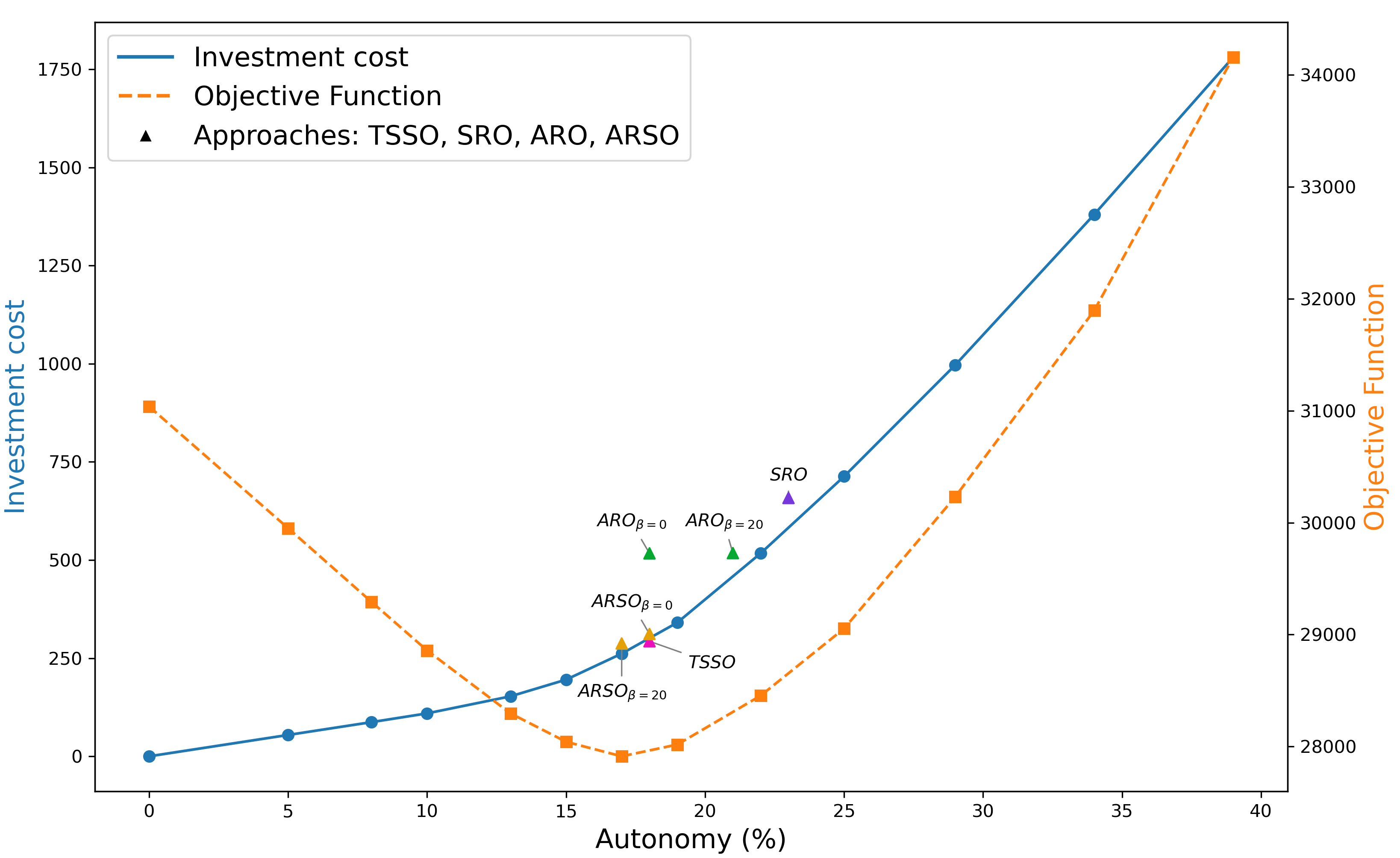}
\caption{Stochastic and robust approaches versus perfect information}
\label{Scatter_Plot}
\end{figure}

To complement our analysis, we have evaluated the corresponding autonomy levels provided by the installed capacities. This enables us to examine where the robust solutions lie on the investment-versus-autonomy curve, offering deeper insight into the trade-offs between investment cost and system autonomy, as shown in Figure~\ref{Scatter_Plot}. To construct this figure, we used the installed capacities of BESS and PV systems determined for each optimization model and ran the operational model over the 216-hour scenario. This allowed us to calculate the energy required from the external grid to meet the electricity demand and estimate the achieved autonomy level for each approach. The triangles in Figure~\ref{Scatter_Plot} represent these results for the TSSO, SRO, ARO, and ARSO formulations, illustrating their respective positions on the investment-versus-autonomy curve.

Additionally, we generated the investment and objective function curves (blue and orange lines in Figure~\ref{Scatter_Plot}) by running a deterministic model under the same 216-hour scenario. For this, we forced specific autonomy levels, ranging from minimal autonomy to full autonomy, and calculated the required DER capacities and associated costs to achieve each level. This approach provides a reference to evaluate how close the robust and stochastic optimization solutions are to the trade-off curve under PI.

Thus, Figure~\ref{Scatter_Plot} shows that ARSO solutions closely approach the optimal solution under the PI scenario, achieving similar levels of autonomy. This indicates a significant reduction in conservatism compared to the ARO solutions. The key to this improvement lies in ARSO's hybrid treatment of uncertainties: addressing ST uncertainty through discrete scenarios while managing LT uncertainty using a BoU approach. This differentiated strategy enables ARSO to balance robustness and cost-efficiency better, bringing it closer to the ideal performance achieved under PI.

\subsection{Algorithm performance}
Table \ref{Convergence_Times} presents the objective function values for the ARO and ARSO formulations and their algorithmic performance in terms of execution time in seconds and total iterations required. Regarding the objective function, when the BoU is zero, the ARO produces a less conservative solution than the ARSO. This behavior is expected because the ARSO incorporates multiple scenarios for PV generation, while the ARO, under a BoU of zero, only considers the predicted value, simplifying the uncertainty representation. However, as $\beta$ increases, the ARSO becomes less conservative than the ARO. This occurs because the ARSO leverages its hybrid nature by combining stochastic scenarios for PV generation with robust methods for demand. 

Concerning the algorithm performances, the ARSO formulation requires more iterations and longer computational times to converge compared to the ARO. This is because the ARSO incorporates both stochastic and robust elements, increasing the problem's complexity. Specifically, the need to evaluate multiple scenarios in the stochastic component adds additional layers to the solution process, particularly as $\beta$ increases. A higher $\beta$ value leads to a broader uncertainty set, generating more dual cuts and thereby increasing the number of iterations and computational effort required for convergence. In contrast, the ARO operates with a more straightforward structure, focusing on demand and generation uncertainties within the BoU framework, resulting in faster convergence. 


\begin{table}[htpb]
\centering
\caption{Stochastic and robust optimization formulations performance}\label{Convergence_Times}
\resizebox{0.6\textwidth}{!}{%
\begin{tabular}{lcrc}
\hline
\textbf{Model}             & \textbf{OF[\euro]} & \textbf{Time[s]} & \textbf{Nº of Iterations} \\ \hline
\textbf{$ARO_{\beta=0}$}   & 2809   & 48         & 23                       \\
\textbf{$ARO_{\beta=10}$}  & 3271     & 141        & 25                       \\
\textbf{$ARO_{\beta=20}$}  & 3519     & 399       & 29                       \\
\textbf{$ARSO_{\beta=0}$}  & 2839    & 364        & 51                       \\
\textbf{$ARSO_{\beta=10}$} & 3193     & 3912       & 57                       \\
\textbf{$ARSO_{\beta=20}$} & 3441    & 4841       & 71                       \\ \hline
\end{tabular}%
}
\end{table}

Table \ref{MP_SP_Times} presents a summary of the computation times in seconds of all iterations for the ARSO and ARO approaches. For each BoU value, the minimum, maximum, and average time of the iterations is presented first, followed by the time of the first and last iteration, respectively. The results clearly show that the algorithm spends significantly more time solving the subproblem, with this time increasing as $\beta$ grows. This increase is particularly pronounced for the ARSO approach, where the subproblem incorporates discrete scenarios, adding complexity to the solution process. In contrast, the computational time for the MP in the ARSO approach remains relatively stable, even as additional Benders-Dual cuts are introduced. This behavior suggests that the computational burden in ARSO is primarily driven by the complexity of the subproblem rather than the iterative updates in the MP. Meanwhile, the ARO approach shows a more balanced increase in both the MP and SP computational times as $\beta$ increases. These observations highlight the importance of properly managing subproblem complexity, especially when discrete scenarios are involved, to improve the scalability and efficiency of the optimization algorithm.

\setlength{\tabcolsep}{4pt}
\begin{table}[htpb]
\centering
\caption{Summary of the computation times in seconds of all iterations for the ARSO and ARO formulation, and BoU values}\label{MP_SP_Times}
\resizebox{0.65\textwidth}{!}{%
\begin{tabular}{clrrrrr}
\hline
\textbf{Problem} & \multicolumn{1}{c}{\textbf{Formulation}}  & \multicolumn{1}{c}{\textbf{Min}} & \multicolumn{1}{c}{\textbf{Max}} & \multicolumn{1}{c}{\textbf{Average}} & \multicolumn{1}{c}{\textbf{First}} & \multicolumn{1}{c}{\textbf{Last}} \\ \hline
MP               & $ARSO_{\beta=0}$          & 2.40          & 7.49         & 4.30             & 2.76             & 5.17           \\
                 & $ARSO_{\beta=10}$          & 2.89         & 7.40          & 4.49             & 2.89             & 6.57           \\
                 & $ARSO_{\beta=20}$          & 2.57         & 7.33         & 4.42             & 2.78             & 6.84           \\
                 & $ARO_{\beta=0}$            & 0.72         & 1.16         & 1.00             & 0.82             & 1.11           \\
                 & $ARO_{\beta=10}$           & 0.76         & 1.20          & 1.03             & 0.78             & 1.15           \\
                 & $ARO_{\beta=20}$           & 0.76         & 1.21         & 1.00             & 0.76             & 1.07           \\ \hline
SP               & $ARSO_{\beta=0}$           & 2.68         & 3.34         & 2.85             & 3.26             & 2.74           \\
                 & $ARSO_{\beta=10}$          & 52.8        & 71.90        & 59.60            & 58.60            & 56.20          \\
                 & $ARSO_{\beta=20}$          & 55.4         & 79.70        & 63.76            & 62.70            & 57.90          \\
                 & $ARO_{\beta=0}$            & 0.89         & 1.30          & 1.07             & 0.91             & 1.07           \\
                 & $ARO_{\beta=10}$          & 3.90          & 5.80          & 4.63             & 5.78             & 5.80            \\
                 & $ARO_{\beta=20}$           & 5.30          & 32.40        & 12.80            & 15.10            & 31.60          \\ \hline
\end{tabular}%
}
\end{table}

\newpage
\section{Conclusions}\label{sec:conclusions}
This study explored and compared different robust and stochastic for the OSP-DER in DNs, considering both LT and ST uncertainties. ARO and ARSO models were proposed and tested on a modified IEEE 33-bus system, where electricity demand was treated as an LT uncertainty and PV generation as an ST uncertainty. The results suggest that the ARSO methodology mitigates the conservatism observed in ARO solutions, achieving a better balance between robustness and cost while more closely approximating the optimal solution under PI. In this regard, the ARSO formulation appears to offer greater flexibility by differentiating between LT demand and ST PV generation uncertainty, which may lead to more accurate investment decisions in DERs and help avoid excessive capacity allocation. Furthermore, the adaptation of the Benders decomposition algorithm to efficiently handle budgeted uncertainty sets proved effective in solving the dual subproblems of ARO and ARSO. This approach successfully managed the problem’s complexity, ensuring computational tractability and achieving convergence within a reasonable number of iterations.

Future research could explore the possibility of relaxing the assumptions of the lossless network model or including market constraints that could influence the investment and operation of DERs. In addition, integrating emerging technologies, such as electric vehicles, or considering other energy carriers, such as hydrogen or thermal energy, could provide a more comprehensive approach to energy planning. Furthermore, a robust data-driven optimization framework could be proposed for the same model, leveraging historical consumption and data generation to build uncertainty sets. This would allow a direct comparison between data-driven approaches and traditional robust optimization approaches, providing insights into the impact of uncertainty representation on solutions. Finally, our Benders decomposition algorithms, such as acceleration techniques or heuristics, require further improvements to deal with large-scale instances.

\section*{Acknowledgment}
This work has been supported by ANID FONDECYT Iniciación 11240745.

\bibliographystyle{APA_Elsevier}
\bibliography{references}


\appendix

\section{Extended constraints of  the ARO subproblem}\label{appendix:spARO}
{\small 
\begin{flalign*} 
    &\displaystyle\sum_{(i,j) \in {\mathcal{L}}} p^{i,j}_{t} - \displaystyle\sum_{(j,i) \in {\mathcal{L}}} p^{j,i}_{t} - pg_{i,t} - pv_{i,t} - ds_{i,t} + ch_{i,t} = -\widetilde{PL}_{i,t} &&: \mathbf{a}_{i,t} \in \mathbb{R} &&& \forall i \in {\mathcal {B}}, \forall t \in {\mathcal {T}}\\
    &\displaystyle\sum_{(i,j) \in {\mathcal{L}}} q^{i,j}_{t} - \displaystyle\sum_{(j,i) \in {\mathcal{L}}} q^{j,i}_{t} - qg_{i,t} = -QL_{i,t} &&: \mathbf{b}_{i,t} \in \mathbb{R} &&& \forall i \in {\mathcal {B}}, \forall t \in {\mathcal {T}}\\
    &v_{j,t} - v_{i,t} + 2R_{i,j} p^{i,j}_{t} + 2X_{i,j}q^{i,j}_{t} = 0  &&: \mathbf{c}_{i,j,t} \in \mathbb{R} &&& \forall (i,j) \in {\mathcal {L}}, \forall t \in {\mathcal {T}} \\  
    &\theta_{i,t} -  \theta_{j,t} - X_{i,j} p^{i,j}_{t} + R_{i,j}q^{i,j}_{t} = 0 &&: \mathbf{d}_{i,j,t} \in \mathbb{R} &&& \forall (i,j) \in {\mathcal {L}}, \forall t \in {\mathcal {T}} \\
    &p^{i,j}_{t}A_{r} + q^{i,j}_{t}B_{r} \leq -S^{max}_{i,j}C_{r} &&: \mathbf{e}^{i,j}_{t,r} \leq 0 &&& \forall (i,j) \in {\mathcal {L}}, \forall t \in {\mathcal {T}}, \forall r \in {\mathcal {R}} \\
    &pg_{i,t} \leq PG^{max}_{i} &&: \mathbf{f}_{i,t} \leq 0  &&& \forall i \in {\mathcal {B}}, \forall t \in {\mathcal {T}}\\
    &qg_{i,t} \geq QG^{min}_{i} &&: \mathbf{g}^{min}_{i,t} \geq 0 &&& \forall i \in {\mathcal {B}}, \forall t \in {\mathcal {T}}\\
    &qg_{i,t} \leq QG^{max}_{i} &&: \mathbf{g}^{max}_{i,t} \leq 0 &&& \forall i \in {\mathcal {B}}, \forall t \in {\mathcal {T}}\\
    &v_{i,t} \geq V^{min}_{i} &&: \mathbf{h}^{min}_{i,t} \geq 0 &&& \forall i \in {\mathcal {B}}, \forall t \in {\mathcal {T}}\\
    &v_{i,t} \leq V^{max}_{i} &&: \mathbf{h}^{max}_{i,t} \leq 0 &&& \forall i \in {\mathcal {B}}, \forall t \in {\mathcal {T}}\\ 
    &pv_{i,t} \leq \widetilde{PV}_{t} \widehat{\gamma}^{pv}_{i}  &&: \mathbf{i}_{i,t} \leq 0 &&& \forall i \in {\mathcal {B}}, \forall t \in {\mathcal {T}} 	\\
	&soc_{i,t} - soc_{i,t-1} - \varphi^{ch} ch_{i,t} + \frac{1}{\varphi^{ds}}ds_{i,t} = 0 &&: \mathbf{j}_{i,t} \in \mathbb{R} &&& \forall i \in {\mathcal {B}}, \forall t \in {\mathcal {T}}, \text{if} \enspace t>1\\
    &soc_{i,1} - \varphi^{ch} ch_{i,1} + \frac{1}{\varphi^{ds}}ds_{i,1} = SOC^{init}\widehat{\gamma}^{bt}_{i} &&: \mathbf{j}_{i,1} \in \mathbb{R} &&& \forall i \in {\mathcal {B}},, \text{if} \enspace t=1 \\
	&soc_{i,t} \geq SOC^{min} \widehat{\gamma}^{bt}_{i} &&: \mathbf{k}^{min}_{i,t} \geq 0 &&& \forall i \in {\mathcal {B}}, \forall t \in {\mathcal {T}}\\
    &soc_{i,t} \leq SOC^{max} \widehat{\gamma}^{bt}_{i} &&: \mathbf{k}^{max}_{i,t} \leq 0 &&& \forall i \in {\mathcal {B}}, \forall t \in {\mathcal {T}}\\
	&ch_{i,t}  \leq PB (\widehat{w}_{i,t}) &&: \mathbf{l}_{i,t} \leq 0 &&& \forall i \in {\mathcal {B}}, \forall t \in {\mathcal {T}} \\ 
	&ds_{i,t} \leq PB (1-\widehat{w}_{i,t}) &&: \mathbf{m}_{i,t} \leq 0 &&& \forall i \in {\mathcal {B}}, \forall t \in {\mathcal {T}}\\
    &\theta_{i,t} \geq \Theta^{min}_{i} &&: \mathbf{n}^{min}_{i,t} \geq 0 &&& \forall i \in {\mathcal {B}}, \forall t \in {\mathcal {T}}\\
    &\theta_{i,t} \leq \Theta^{max}_{i} &&: \mathbf{n}^{max}_{i,t} \leq 0 &&& \forall i \in {\mathcal {B}}, \forall t \in {\mathcal {T}}
\end{flalign*}
}

\newpage

\section{Extended formulation of the ARO dual subproblem}\label{appendix:dspARO}

{\small 
\begin{multline}
\nonumber \underset{\lambda, \mu }{max} \displaystyle\sum_{i\in {\mathcal {B}}} \displaystyle\sum_{t\in {\mathcal {T}}} \bigg\{   -(\overline{PL}_{i,t}\mathbf{a}_{i,t} + \widehat{PL}_{i,t}(\mathbf{a}^{+}_{i,t} - \mathbf{a}^{-}_{i,t})) - QL_{i,t}\mathbf{b}_{i,t} + 
PG^{max}_{i}\mathbf{f}_{i,t} + 
\mathbf{g}^{min}_{i,t}QG^{min}_{i} +  \mathbf{g}^{max}_{i,t}QG^{max}_{i} \\ 
+ \mathbf{h}^{min}_{i,t}V^{min}_{i}  +  \mathbf{h}^{max}_{i,t}V^{max}_{i} + 
(\overline{PV}_{t} + \widehat{PV}_{t}(\mathbf{i}^{+}_{i,t} - \mathbf{i}^{-}_{i,t})) \widehat{\gamma}^{pv}_{i} + \mathbf{k}^{min}_{i,t}SOC^{min}\widehat{\gamma}^{bt}_{i} 
+  \mathbf{k}^{max}_{i,t}SOC^{max} \widehat{\gamma}^{bt}_{i} \\ + 
\mathbf{l}_{i,t}PB\widehat{w}_{i,t}   + \mathbf{m}_{i,t}(PB(1-\widehat{w}_{i,t})) +  \mathbf{n}^{min}_{i,t}\Theta^{min}_{i} + 
\mathbf{n}^{max}_{i,t}\Theta^{max}_{i}     \bigg\} \\ 
- \displaystyle\sum_{(i,j)\in {\mathcal {L}}} \displaystyle\sum_{r\in {\mathcal {R}}} S^{max}_{i,j}C_{r}\mathbf{e}^{r}_{i,j,t} + \displaystyle\sum_{i\in {\mathcal {B}}}  SOC^{init} \widehat{\gamma}^{bt}_{i} \mathbf{j}_{i,1} 
\end{multline}

\begin{flalign*} 
    &\mathbf{f}_{i,t} -\mathbf{a}_{i,t} \leq \rho_s \lambda_{t} && \forall i \in {\mathcal {B}}, \forall t \in {\mathcal {T}}\\
    &\mathbf{i}_{i,t} -\mathbf{a}_{i,t} \leq \rho_s OC^{pv} && \forall i \in {\mathcal {B}}, \forall t \in {\mathcal {T}}\\
    &\frac{\mathbf{j}_{i,t}}{\varphi} -\mathbf{a}_{i,t} + \mathbf{m}_{i,t} \leq \rho_s OC^{bt} && \forall i \in {\mathcal {B}}, \forall t \in {\mathcal {T}}\\
    &\mathbf{a}_{i,t} - \varphi\mathbf{j}_{i,t} + \mathbf{l}_{i,t} \leq 0&& \forall i \in {\mathcal {B}}, \forall t \in {\mathcal {T}} \\
    &\mathbf{j}_{i,t} - \mathbf{j}_{i,t+1} + \mathbf{k}^{min}_{i,t} + \mathbf{k}^{max}_{i,t}  \leq 0 && \forall i \in {\mathcal {B}}, \forall t \in {\mathcal {T}}\\
    &\mathbf{g}^{min}_{i,t} + \mathbf{g}^{max}_{i,t} - \mathbf{b}_{i,t} = 0 && \forall i \in {\mathcal {B}}, \forall t \in {\mathcal {T}}\\
    &\mathbf{h}^{min}_{i,t} + \mathbf{h}^{max}_{i,t} + \displaystyle\sum_{(i,j)\in {\mathcal {L}}} (\mathbf{c}_{i,j,t} - \mathbf{c}_{j,i,t}) \leq 0 && \forall i \in {\mathcal {B}}, \forall t \in {\mathcal {T}}\\
    &\mathbf{n}^{min}_{i,t} + \mathbf{n}^{max}_{i,t} + \displaystyle\sum_{(i,j)\in {\mathcal {L}}} (\mathbf{d}_{j,i,t} - \mathbf{d}_{i,j,t}) = 0 && \forall i \in {\mathcal {B}}, \forall t \in {\mathcal {T}}\\
    &\mathbf{a}_{i,t} - \mathbf{a}_{j,t} + 2R_{i,j}\mathbf{c}_{i,j,t} - X_{i,j}\mathbf{d}_{i,j,t}  + \displaystyle\sum_{r\in {\mathcal {R}}} (A_{r}\mathbf{e}^{r}_{i,j,t}) = 0 && (i,j)\in {\mathcal {L}}, \forall t \in {\mathcal {T}}\\
    &\mathbf{b}_{i,t} - \mathbf{b}_{j,t} + 2X_{i,j}\mathbf{c}_{i,j,t} + R_{i,j}\mathbf{d}_{i,j,t} + \displaystyle\sum_{r\in {\mathcal {R}}}(B_{r}\mathbf{e}^{r}_{i,j,t}) = 0 && (i,j)\in {\mathcal {L}}, \forall t \in {\mathcal {T}} \\
    &\displaystyle\sum_{t\in {\mathcal {T}}} (U^{+}_{t} + U^{-}_{t}) \leq \beta_{pv} \\
    &\displaystyle\sum_{t\in {\mathcal {T}}} (V^{+}_{i,t} + V^{-}_{i,t}) \leq \beta_{pl} &&\forall i \in {\mathcal {B}}\\
    &U^{+}_{t} + U^{-}_{t} \leq 1 && \forall t \in {\mathcal {T}}\\
    &V^{+}_{i,t} + V^{-}_{i,t} \leq 1 && \forall i \in {\mathcal {B}}, \forall t \in {\mathcal {T}}\\
    &-V^{+}_{i,t} M \leq \mathbf{a}^{+}_{i,t} \leq V^{+}_{i,t} M  && \forall i \in {\mathcal {B}}, \forall t \in {\mathcal {T}}\\
    &-V^{-}_{i,t} M \leq \mathbf{a}^{-}_{i,t} \leq V^{-}_{i,t} M && \forall i \in {\mathcal {B}}, \forall t \in {\mathcal {T}} \\
    &-U^{+}_{t} M \leq \mathbf{i}^{+}_{i,t} \leq U^{+}_{t} M  && \forall i \in {\mathcal {B}}, \forall t \in {\mathcal {T}}\\
    &-U^{-}_{t} M \leq \mathbf{i}^{-}_{i,t} \leq U^{-}_{t} M  && \forall i \in {\mathcal {B}}, \forall t \in {\mathcal {T}}\\
    &\mathbf{a}_{i,t} - M(1-V^{+}_{i,t}) \leq \mathbf{a}^{+}_{i,t} \leq \mathbf{a}_{i,t} + M(1-V^{+}_{i,t})  && \forall i \in {\mathcal {B}}, \forall t \in {\mathcal {T}}\\
    &\mathbf{a}_{i,t} - M(1-V^{-}_{i,t}) \leq \mathbf{a}^{-}_{i,t} \leq \mathbf{a}_{i,t} + M(1-V^{-}_{i,t})  && \forall i \in {\mathcal {B}}, \forall t \in {\mathcal {T}}\\
    &\mathbf{i}_{i,t} - M(1-U^{+}_{t}) \leq \mathbf{i}^{+}_{i,t} \leq \mathbf{i}_{i,t} + M(1-U^{+}_{t})  && \forall i \in {\mathcal {B}}, \forall t \in {\mathcal {T}}\\
    &\mathbf{i}_{i,t} - M(1-U^{-}_{t}) \leq \mathbf{i}^{-}_{i,t} \leq \mathbf{i}_{i,t} + M(1-U^{-}_{t})  && \forall i \in {\mathcal {B}}, \forall t \in {\mathcal {T}}
\end{flalign*}
}

\newpage

\section{Extended formulation of the ARSO dual subproblem}\label{appendix:dspARSO}

\begin{multline}
\nonumber \underset{\lambda, \mu }{max} \displaystyle\sum_{i\in {\mathcal {B}}} \displaystyle\sum_{t\in {\mathcal {T}}} \displaystyle\sum_{s\in {\mathcal {S}}}\bigg\{   -(\overline{PL}_{i,t}\mathbf{a}_{i,t,s} + \widehat{PL}_{i,t}(\mathbf{a}^{+}_{i,t,s} - \mathbf{a}^{-}_{i,t,s})) - QL_{i,t,s}\mathbf{b}_{i,t,s} + 
PG^{max}_{i}\mathbf{f}_{i,t,s} + 
\mathbf{g}^{min}_{i,t,s}QG^{min}_{i} \\ +  \mathbf{g}^{max}_{i,t,s}QG^{max}_{i}  
+ \mathbf{h}^{min}_{i,t,s}V^{min}_{i} +  \mathbf{h}^{max}_{i,t,s}V^{max}_{i} + 
PV_{t,s} \widehat{\gamma}^{pv}_{i}\mathbf{i}_{i,t,s} + \mathbf{k}^{min}_{i,t,s}SOC^{min}\widehat{\gamma}^{bt}_{i} 
+  \mathbf{k}^{max}_{i,t,s}SOC^{max} \widehat{\gamma}^{bt}_{i} \\ + 
\mathbf{l}_{i,t,s}PB\widehat{w}_{i,t,s} + \mathbf{m}_{i,t,s}(PB(1-\widehat{w}_{i,t,s})) +  \mathbf{n}^{min}_{i,t,s}\Theta^{min}_{i} + 
\mathbf{n}^{max}_{i,t,s}\Theta^{max}_{i}     \bigg\} \\ 
- \displaystyle\sum_{(i,j)\in {\mathcal {L}}} \displaystyle\sum_{r\in {\mathcal {R}}} \displaystyle\sum_{s\in {\mathcal {S}}} S^{max}_{i,j}C_{r}\mathbf{e}^{r}_{i,j,t,s} + \displaystyle\sum_{i\in {\mathcal {B}}} \displaystyle\sum_{s\in {\mathcal {S}}} SOC^{init} \widehat{\gamma}^{bt}_{i} \mathbf{j}_{i,1,s} 
\end{multline}

\begin{flalign*} 
    &\mathbf{f}_{i,t,s} -\mathbf{a}_{i,t,s} \leq \rho_s \lambda_{t} && \forall i \in {\mathcal {B}}, \forall t \in {\mathcal {T}}, \forall s \in {\mathcal {S}}\\
    &\mathbf{i}_{i,t,s} -\mathbf{a}_{i,t,s} \leq \rho_s OC^{pv} && \forall i \in {\mathcal {B}}, \forall t \in {\mathcal {T}}, \forall s \in {\mathcal {S}} \\
    &\frac{\mathbf{j}_{i,t,s}}{\varphi} -\mathbf{a}_{i,t,s} + \mathbf{m}_{i,t,s} \leq \rho_s OC^{bt}  && \forall i \in {\mathcal {B}}, \forall t \in {\mathcal {T}}, \forall s \in {\mathcal {S}}\\
    &\mathbf{a}_{i,t,s} - \varphi\mathbf{j}_{i,t,s} + \mathbf{l}_{i,t,s} \leq 0 && \forall i \in {\mathcal {B}}, \forall t \in {\mathcal {T}}, \forall s \in {\mathcal {S}}\\
    &\mathbf{j}_{i,t,s} - \mathbf{j}_{i,t+1,s} + \mathbf{k}^{min}_{i,t,s} + \mathbf{k}^{max}_{i,t,s}  \leq 0 && \forall i \in {\mathcal {B}}, \forall t \in {\mathcal {T}}, \forall s \in {\mathcal {S}}\\
    &\mathbf{g}^{min}_{i,t,s} + \mathbf{g}^{max}_{i,t,s} - \mathbf{b}_{i,t,s} = 0 && \forall i \in {\mathcal {B}}, \forall t \in {\mathcal {T}}, \forall s \in {\mathcal {S}}\\
    &\mathbf{h}^{min}_{i,t,s} + \mathbf{h}^{max}_{i,t,s} + \displaystyle\sum_{(i,j)\in {\mathcal {L}}} (\mathbf{c}_{i,j,t,s} - \mathbf{c}_{j,i,t,s}) \leq 0 && \forall i \in {\mathcal {B}}, \forall t \in {\mathcal {T}}, \forall s \in {\mathcal {S}}\\
    &\mathbf{n}^{min}_{i,t,s} + \mathbf{n}^{max}_{i,t,s} + \displaystyle\sum_{(i,j)\in {\mathcal {L}}} (\mathbf{d}_{j,i,t,s} - \mathbf{d}_{i,j,t,s}) = 0&& \forall i \in {\mathcal {B}}, \forall t \in {\mathcal {T}}, \forall s \in {\mathcal {S}}\\
    &\mathbf{a}_{i,t,s} - \mathbf{a}_{j,t,s} + 2R_{i,j}\mathbf{c}_{i,j,t,s} - X_{i,j}\mathbf{d}_{i,j,t,s}  + \displaystyle\sum_{r\in {\mathcal {R}}} (A_{r}\mathbf{e}^{r}_{i,j,t,s}) = 0&& (i,j)\in {\mathcal {L}}, \forall t \in {\mathcal {T}}, \forall s \in {\mathcal {S}}\\
    &\mathbf{b}_{i,t,s} - \mathbf{b}_{j,t,s} + 2X_{i,j}\mathbf{c}_{i,j,t,s} + R_{i,j}\mathbf{d}_{i,j,t,s} + \displaystyle\sum_{r\in {\mathcal {R}}}(B_{r}\mathbf{e}^{r}_{i,j,t,s}) = 0 && (i,j)\in {\mathcal {L}}, \forall t \in {\mathcal {T}}, \forall s \in {\mathcal {S}}\\
    &\displaystyle\sum_{t\in {\mathcal {T}}} (V^{+}_{i,t} + V^{-}_{i,t}) \leq \beta_{pl} &&\forall i \in {\mathcal {B}}\\
    &V^{+}_{i,t} + V^{-}_{i,t} \leq 1 && \forall i \in {\mathcal {B}}, \forall t \in {\mathcal {T}}\\
    &-V^{+}_{i,t} M \leq \mathbf{a}^{+}_{i,t,s} \leq V^{+}_{i,t} M  && \forall i \in {\mathcal {B}}, \forall t \in {\mathcal {T}}, \forall s \in {\mathcal {S}}\\
    &-V^{-}_{i,t} M \leq \mathbf{a}^{-}_{i,t,s} \leq V^{-}_{i,t} M && \forall i \in {\mathcal {B}}, \forall t \in {\mathcal {T}}, \forall s \in {\mathcal {S}} \\
    &\mathbf{a}_{i,t,s} - M(1-V^{+}_{i,t}) \leq \mathbf{a}^{+}_{i,t,s} \leq \mathbf{a}_{i,t,s} + M(1-V^{+}_{i,t})  && \forall i \in {\mathcal {B}}, \forall t \in {\mathcal {T}}, \forall s \in {\mathcal {S}}\\
    &\mathbf{a}_{i,t,s} - M(1-V^{-}_{i,t}) \leq \mathbf{a}^{-}_{i,t,s} \leq \mathbf{a}_{i,t,s} + M(1-V^{-}_{i,t})  && \forall i \in {\mathcal {B}}, \forall t \in {\mathcal {T}}, \forall s \in {\mathcal {S}}
\end{flalign*}

\end{document}